\documentclass[11pt]{article}
\usepackage{amscd,amsmath,amssymb}
\usepackage[mathscr]{eucal}

\makeatletter \@addtoreset{equation}{section}\makeatother

\title{\bf Covariant $q$-differential operators and unitary
highest weight representations  for $U_q\mathfrak{su}_{n,n}$}
\author{Dmitry Shklyarov$^\dagger$ \and Genkai Zhang$^\ddagger$}
\date{}

\begin{document}
\maketitle

\centerline{$^\dagger$\small{Department of Mathematics, Kansas State University}} \centerline {\small{Manhattan, KS 66506, USA}}

\centerline{$^\ddagger$\small{Chalmers Tekniska H\"ogskola/G\"oteborgs
Universitet, Matematik}} \centerline{\small{412 96, G\"oteborg, Sweden}}

\begin{center}
\small{e-mail: shklyarov@math.ksu.edu, genkai@math.chalmers.se}
\end{center}

\makeatletter
\let\@thefnmark\relax
\makeatother

\newtheorem{theorem}{Theorem}[section]
\newtheorem{lemma}[theorem]{Lemma}
\newtheorem{proposition}[theorem]{Proposition}
\newtheorem{corollary}[theorem]{Corollary}
\bigskip
\begin{center}\begin{minipage}[t]{4in}\small{ABSTRACT:
 We investigate a
one-parameter family of quantum Harish-Chandra modules of
$U_q\mathfrak{sl}_{2n}$. This family is 
an analog of the holomorphic discrete series
of representations of the group $SU(n,n)$
for the quantum group $U_q\mathfrak{su}_{n, n}$.  We introduce a $q$-analog of "the
wave" operator (a determinant-type differential operator) and prove certain
covariance property of its powers. This result is applied to the study of
some quotients of the above-mentioned quantum Harish-Chandra modules. We also prove
an analog of a known result by J.~Faraut and A.~Koranyi 
on the expansion of reproducing kernels which 
determines the analytic continuation of the holomorphic discrete series.}
\end{minipage}
\end{center}

\bigskip
\begin{center}\begin{minipage}[t]{4in}\small{\centerline{Table of contents}

1. Introduction

2. Quantum space of matrices and its symmetries

$\hspace{5mm}$ 2.1 Quantum space of matrices

$\hspace{5mm}$ 2.2 A structure of $U_q\mathfrak{sl}_{2n}$-module algebra on
$\mathbb{C}[M_{n}]_q$

$\hspace{5mm}$ 2.3 A twisted $U_q\mathfrak{sl}_{2n}$-action on
$\mathbb{C}[M_{n}]_q$

3. Some $q$-differential operators

$\hspace{5mm}$ 3.1 Basic definitions

$\hspace{5mm}$ 3.2 A $q$-wave operator

4. A covariance property

$\hspace{5mm}$ 4.1 Formulation

$\hspace{5mm}$ 4.2 A $q$-analog of the Cauchy-Szeg\"o integral
representation

$\hspace{5mm}$ 4.3  $q$-Analogs of the kernels $\mathrm{det}
(1-\mathbf{z}\boldsymbol{\zeta}^*)^{-N}$

$\hspace{5mm}$ 4.4 A $q$-analog of the Fock inner product

$\hspace{5mm}$ 4.5 Proof of the covariance property

5. Holomorphic discrete series for $U_q\mathfrak{su}_{n,n}$

$\hspace{5mm}$ 5.1 Definition of the holomorphic discrete series

$\hspace{5mm}$ 5.2  A $q$-analog of a result by J.~Faraut and A.~Koranyi

$\hspace{5mm}$ 5.3 Some consequences of the covariance property}

References
\end{minipage}
\end{center}

\bigskip

\section{Introduction}

We start with recalling classical known results about analytic continuation
of the weighted Bergman spaces in the unit disk and their explicit
realization.

Recall that the group $SU_{1,1}$ acts on the unit disk by fractional-linear
transformations. Many important representations of the group are realized
geometrically in various functional spaces on the disk and on the unit
circle. In particular, representations of the discrete series admit a
realization of that kind. Namely, consider the kernel
$(1-z\overline{w})^{-\lambda}$ in the unit disk. For $\lambda>1$ it is the
reproducing kernel for the so-called weighted Bergman space consisting of
holomorphic functions that are square integrable with the weight
$(1-|z|^2)^{\lambda-2}dm(z)$ (here $dm(z)$ is the normalized Lebesgue
measure). The group $SU_{1,1}$ acts in the space via change of variable and
a multiplier:
\begin{equation}\label{tw}
\pi_\lambda(g)(f(z))=f(g^{-1}z)\cdot(cz+d)^{-\lambda},\quad g^{-1}=
\left(\begin{array}{cl} a & \quad b\\ c & \quad d
\end{array}\right)
\end{equation}
(for non-integer $\lambda$'s one should consider the universal covering
$\widetilde{SU}_{1,1}$ instead of $SU_{1,1}$). Thus obtained unitary representation
belongs to the discrete series and is said to be a representation of the
holomorphic discrete series for $SU_{1,1}$ or $\widetilde{SU}_{1,1}$.

The reproducing kernel $(1-z\overline{w})^{-\lambda}$ has analytic
continuation in the parameter $\lambda$. This is obtained from the formula
\begin{equation}\label{de}
(1-z\overline{w})^{-\lambda}=\sum\limits_{m=0}^{\infty}(\lambda)_m
\frac{(z\overline{w})^m}{m!},\quad
(\lambda)_m=\lambda\cdot(\lambda+1)\cdot\ldots\cdot(\lambda+m-1).
\end{equation}

For $\lambda>0$ the kernel is still positive definite, and the
$SU_{1,1}$-action (\ref{tw}) in the associated Hilbert space is also
unitary. For $\lambda=1$, the Hilbert space is the Hardy space of
holomorphic function on the closed disk whose boundary value are square
integrable on the circle.

For further study of the previous representations, it is convenient to pass to
the corresponding Harish-Chandra modules. Consider the space $\mathbb{C}[z]$
of polynomials on $\mathbb{C}$. The representation $\pi_\lambda$ induces a
representation of
$\mathfrak{sl}_2(\mathbb{C})=\mathfrak{su}_{1,1}\otimes\mathbb{C}$ on
$\mathbb{C}[z]$ which may be defined for any $\lambda\in\mathbb{R}$ (and
even for $\lambda\in\mathbb{C}$). Let us denote by $\mathcal{P}_{\lambda}$
the space $\mathbb{C}[z]$ endowed with the above-mentioned action of
$\mathfrak{sl}_2(\mathbb{C})$. $\mathcal{P}_{\lambda}$ is irreducible for
all positive $\lambda$'s. However, if $\lambda=1-l$ for some positive integer
$l$ then $\mathcal{P}_{\lambda}=\mathcal{P}_{1-l}$ has the following composition
series
\begin{equation}\label{compo}
\{0\}\subset\mathcal{P}^{(0)}_{1-l}\subset\mathcal{P}_{1-l}
\end{equation}
where $\mathcal{P}^{(0)}_{1-l}$ is the submodule of polynomials of
degree $\leq l-1$. The natural problem is to study the quotient
$\mathcal{P}_{1-l}/\mathcal{P}^{(0)}_{1-l}$. This is the point where
covariant differential operators appear on the scene. They play an important
role in an explicit realization of the quotient. Namely, one checks that
the differential operator $\left(\frac{\partial}{\partial z}\right)^l$
intertwines the actions $\pi_{1-l}$ and $\pi_{1+l}$:
\begin{equation}\label{bolslemma}
\left(\frac{\partial}{\partial z}\right)^l\cdot\pi_{1-l}(\xi)
=\pi_{1+l}(\xi)\cdot\left(\frac{\partial}{\partial z}\right)^l, \quad
\xi\in\mathfrak{sl}_2(\mathbb{C}).
\end{equation}
Clearly, $\left(\frac{\partial}{\partial z}\right)^l$ induces an isomorphism
from $\mathcal{P}_{1-l}/\mathcal{P}^{(0)}_{1-l}$ into $\mathcal{P}_{1+l}$,
and this, in particular, proves unitarizability of the former module.

The unit disk is the simplest example of a bounded symmetric domain
\cite{Arazy}. The above-mentioned results admit appropriate generalization for any
such domain (of course, the group $SU_{1,1}$ is replaced by the group of
biholomorphic automorphisms of the domain under consideration).

For the so-called tube domains, some generalizations of the covariance property
(\ref{bolslemma}) have been obtained by G.~Shimura \cite{Sh}, J.~Arazy
\cite{Arazy}, H.P.~Jakobsen \cite{Jac}, H.P.~Jakobsen and M.~Vergne
\cite{JacVer}, H.P.~Jakobsen and M.~Harris \cite{JacHar}. For example, in
the case of the tube domain of type $\mathrm{I}_{n,n}$ (the unit
ball in the space of complex $n\times n$-matrices) the analog of (\ref{bolslemma}) is a statement about an intertwining
property of powers of the operator $\Box=\det\left(\frac{\partial}{\partial
z_a^\alpha}\right)_{a,\alpha=1,\ldots n}$ with respect to certain "twisted" action of the group
$SU_{n,n}$ analogous to (\ref{tw}).

The generalized covariance property (\ref{bolslemma}) has turned out to be useful beyond the problems we
mentioned previously. It has been applied also to computing the Harish-Chandra
homomorphism of invariant differential operators \cite{Genkai}.

Now for symmetric bounded domains the expansion (\ref{de}) has been
found by \O{}rsted \cite{Orsted-tams-80}
for type I matrix domains and in general case by
J.~Faraut and A.~Koranyi  \cite{FK}. From this expansion one can
read off the composition series analogous to (\ref{compo});
the covariant property of the intertwining operators
is related to the classical Cayley-Capelli type formula. We
note that the unitarity of the highest weight modules
had been classified earlier by Jakobsen \cite{Jak-clf} using
algebraic method; however the analytic approach as
in \cite{Orsted-tams-80}  and \cite{FK} generated
some other interesting analytic subjects and is related
to many problems in special functions and orthogonal
polynomials. For quantum groups
the classification of unitary highest weight
representations has also been done recently \cite{Jac-letter-mp},
and we believe however that an analytic and concrete approach 
deserves pursuing.

In the present paper we obtain analogs of (\ref{de}), (\ref{compo}), and
(\ref{bolslemma}) for a quantum matrix ball, an analog of the tube domain of
type $\mathrm{I}_{n,n}$ which has been defined in framework of quantum group
theory by L.~Vaksman et al \cite{SSV1}.

In \cite{SSV2}, the authors defined analogs of the weighted Bergman spaces
on the quantum matrix ball. Also, they constructed analogs of the
corresponding reproducing kernels and the "twisted" unitary action of the
group $SU_{n,n}$. From the representation theoretic point of view, the paper
\cite{SSV2} presents a $q$-analog of the holomorphic discrete series of the
group $SU_{n,n}$ (more precisely, analogs of the associated Harish-Chandra
modules).

The natural problem now is to investigate those representations,
particularly, to define their "analytic continuation" and to study
composition series of the resulting Harish-Chandra modules. In the case
$n=2$, these problems were treated in \cite{GeomReal}. In the present paper,
we deal with the case of arbitrary $n$.

The role of covariant differential operators in the classical theory of
bounded symmetric domains and related Harish-Chandra modules is very well
known \cite{Jac, JacHar, JacVer, Sh}. Our intention is to bring
covariant $q$-differential operators into the study of quantum
Harish-Chandra modules and thus to demonstrate their importance in the
quantum setting as well\footnote{Note that similar questions have been
already treated in the literature (see \cite{Dob} and, especially,
\cite{Jacdiff}).}. We introduce a determinant-type $q$-differential
operator similar to $\Box$ and prove a $q$-analog of the covariance property. In
the last section, this result is applied to investigation of certain
quotients of the above quantum Harish-Chandra modules.

Another goal of the paper is to obtain an analog of the aforementioned
result by J.~Faraut and A.~Koranyi which has allowed them to solve the
problem of analytic continuation of the holomorphic discrete series in the
classical setting.

As we already mentioned, there is a complete classification of
unitarizable highest-weight modules over quantum groups (see
\cite{Jac-letter-mp}). Thus, neither the holomorphic discrete series of the
quantum group $SU(n,n)$, constructed in \cite{SSV2}, nor its analytic
continuation, obtained in the present paper, give us a new family of
unitary modules. The principal aim of the present paper, as well as papers
\cite{SSV2}, \cite{GeomReal}, is to develop an "analytic and geometric" framework for
studying quantum Harish-Chandra modules related to quantum Cartan domains,
particularly, to show that there are substantial generalizations of known
classical constructions and results connected with the holomorphic discrete
series.

The paper is organized as follows.
Sections 2 and 3 contain some preliminary material. In Section 2, we recall
some basic notions and results of quantum group theory (particularly, the
notion of quantum space of matrices and of the quantized universal
enveloping algebra $U_q\mathfrak{sl}_n$). This is done mainly for the
purpose to set the notation we use further. Also, we recall certain hidden
quantum $U_q\mathfrak{sl}_{2n}$-symmetry of the quantum matrix space
discovered in \cite{SV1}. \footnote{This hidden symmetry was one of the first hints
that there should be a substantial theory of $q$-bounded symmetric domains.
These objects were invented a little later in \cite{SV}.} In the end of
Section 2 we describe a twisted action (depending on a parameter $\lambda$) of $U_q\mathfrak{sl}_{2n}$ on the
quantum matrix space. For 
$\lambda$ large enough the corresponding Harish-Chandra modules are unitarizable
representations of $U_q\mathfrak{su}_{n,n}$, which we call the holomorphic
discrete series due to the previous motivation. Section 3 is devoted to
$q$-differential operators. We recall there the notion of a $q$-differential
operator with constant coefficients and describe certain properties
of the algebra of such operators. Also, we introduce an analog of the
operator $\Box$ and derive its "obvious" quantum symmetry
which amounts to an intertwining property
of the operator with an action of the Hopf
subalgebra $U_q\mathfrak{sl}_{n}\otimes U_q\mathfrak{sl}_{n}\subset
U_q\mathfrak{sl}_{2n}$. This obvious symmetry is extended to a large hidden
symmetry, namely, the intertwining
property of the operator (and of its powers) with the
twisted $U_q\mathfrak{sl}_{2n}$-actions. This covariance property is formulated and proved in Section 4. In
the course of the proof, we use a number of results from the theory of
quantum bounded symmetric domains, in particular, those obtained in
\cite{Gauss} and, especially, results of \cite{CS}. To keep the size of the paper reasonable, we have to be more sketchy in this part of the paper. We omit proofs of those results giving appropriate references instead.  In the last section of
the paper we investigate the holomorphic
discrete series for $U_q\mathfrak{su}_{n,n}$. First of all, we use
computations of Section 4 to produce an analog of the result by J.~Faraut
and A.~Koranyi we mentioned earlier. Then we derive some applications of the
covariance property.

{\bf Acknowledgments.} This research was supported by Royal Swedish Academy of
Sciences under the program "Cooperation between researchers in Sweden and
the former Soviet Union". The authors are indebted to Leonid Vaksman for sharing 
with us many of his results and
ideas. Moreover, results of subsection \ref{FaKo} are joint with him, and we
are grateful to him for generously allowing us to publish those results
here.

\bigskip

\section{Quantum space of matrices and its symmetries}

In this paper, the parameter $q$ is supposed to be a number from the
interval $(0,1)$.

\medskip

\subsection{Quantum space of matrices}\label{qsom}

Let us start with the definition of the algebra $\mathbb{C}[M_{n}]_q$ of
polynomials on the quantum matrix space. It is the unital algebra given by
its generators $z_a^\alpha$ (here $a,\alpha=1,\ldots n$, $a$ is the column
index
 and $\alpha$ is the  row index) and the following relations
\begin{equation}\label{Mn}
z_a^\alpha z_b^\beta=\left \{\begin{array}{ccl}qz_b^\beta z_a^\alpha &,&
a=b\&\;\alpha<\beta \quad{\rm or}\quad a<b \;\&\;\alpha=\beta \\ z_b^\beta
z_a^\alpha &,& a<b \;\&\;\alpha>\beta \\ z_b^\beta
z_a^\alpha+(q-q^{-1})z_a^\beta z_b^\alpha &,& a<b \;\&\;\alpha<\beta
\end{array} \right..
\end{equation}

These commutation relations, along with the relation
\begin{equation}\label{sln}
\mathrm{det}_q(\mathbf{z})=\sum_{s \in
S_n}(-q)^{l(s)}z_{a_1}^{\alpha_{s(1)}}z_{a_2}^{\alpha_{s(2)}}\cdots
z_{a_n}^{\alpha_{s(n)}}=1,
\end{equation}
appeared for the first time in \cite{Drin} as the relations between
generators in the algebra $\mathbb{C}[SL_n]_q$ of regular functions on the
quantum $SL_n$. It was suggested in \cite{FRT} to discard (\ref{sln}) from
the list of relations and to regard (\ref{Mn}) as the defining relations of
the algebra of polynomials on the quantum space of matrices. The algebra
$\mathbb{C}[SL_n]_q$ is then the quotient of $\mathbb{C}[M_{n}]_q$ by the
two-sided ideal generated by the element $\det \nolimits_q (\mathbf{z})-1$
(note that the $q$-determinant $\det \nolimits_q (\mathbf{z})$ belongs to
the center of $\mathbb{C}[M_n]_q$ (\cite[Section 7.3.B]{ChP})). Also, the
algebra $\mathbb{C}[M_n]_q$ is used to define the algebra of regular
functions on the quantum $GL_n$. The latter is just the localization of the
former with respect to the multiplicative system $\det \nolimits_q
(\mathbf{z})^m$, $m=1,2,\ldots$.

The crucial observation concerning the algebra $\mathbb{C}[M_{n}]_q$ was the
discovery of the comultiplication
$$\mathbb{C}[M_{n}]_q\to
\mathbb{C}[M_{n}]_q\otimes\mathbb{C}[M_{n}]_q, \qquad
z_a^\alpha\mapsto\sum_j z_a^j\otimes z_j^\alpha$$ which, along with the
initial multiplication, makes $\mathbb{C}[M_{n}]_q$ into a bialgebra. The
comultiplication maps the $q$-determinant $\mathrm{det}_q(\mathbf{z})$ to
$\mathrm{det}_q(\mathbf{z})\otimes\mathrm{det}_q(\mathbf{z})$ and thus
induces a comultiplication on the algebra $\mathbb{C}[SL_n]_q$. The latter,
along with certain antipode and counit, makes $\mathbb{C}[SL_n]_q$ into a
Hopf algebra.

All the above structures allow one to produce $q$-analogs of the left and
right actions $$ L(g): f(\mathbf{z})\mapsto f(g^{-1}\cdot\mathbf{z}), \qquad
R(g): f(\mathbf{z})\mapsto f(\mathbf{z}\cdot g)$$ of $SL_n$ in
$\mathbb{C}[M_{n}]$. These $q$-analogs are usually described in terms of
comodule algebras \cite{ChP}. However, it is more convenient for us to use
an "infinitesimal" version of those actions which is based on the notion of
the quantum universal enveloping algebra $U_q\mathfrak{sl}_n$ due to
Drinfeld \cite{Drin} and Jimbo \cite{Jimbo}.

First, we recall the definition of $U_q\mathfrak{sl}_n$ (we follow the
 notation of \cite{J}). The quantum universal enveloping
algebra $U_q\mathfrak{sl}_n$ is the unital algebra generated by the elements
$E_i$, $F_i$, $K_i^{\pm1}$, $i=1,\ldots,n$, which satisfy the relations $$
K_iK_j=K_jK_i,\quad K_iK_i^{-1}=K_i^{-1}K_i=1,$$ $$ K_iE_j=q^{a_{ij}}E_jK_i,
\quad K_iF_j=q^{-a_{ij}}F_jK_i, $$ $$
E_iF_j-F_jE_i=\delta_{ij}(K_i-K_i^{-1})/(q-q^{-1}),$$
$$ E_i^2E_j-(q+q^{-1})E_iE_jE_i+E_jE_i^2=0,\quad |i-j|=1 $$ $$
F_i^2F_j-(q+q^{-1})F_iF_jF_i+F_jF_i^2=0,\quad |i-j|=1 $$ $$
[E_i,E_j]=[F_i,F_j]=0,\quad |i-j|\ne 1 $$ with $(a_{ij})$ being the Cartan
matrix of
type $A_{n-1}$. Moreover, $U_q\mathfrak{sl}_n$ is a Hopf algebra. The
comultiplication
$\Delta$,  the antipode $S$, and the counit $\varepsilon$ are determined by
\begin{equation}\label{delta}
\Delta(E_i)=E_i \otimes 1+K_i \otimes E_i,\quad \Delta(F_i)=F_i \otimes
K_i^{-1}+1 \otimes F_i,\quad \Delta(K_i)=K_i \otimes K_i,
\end{equation}
\begin{equation}\label{S}
S(E_i)=-K_i^{-1}E_i,\quad S(F_i)=-F_iK_i,\quad S(K_i)=K_i^{-1},
\end{equation}
\begin{equation}\label{eps}
\varepsilon(E_i)=\varepsilon(F_i)=0,\quad \varepsilon(K_i)=1.
\end{equation}

It is observed in \cite{Drin} that the Hopf algebras $U_q\mathfrak{sl}_n$
and $\mathbb{C}[SL_n]_q$ are dual to each other. This, in particular, allows
one to use the language of $U_q\mathfrak{sl}_n$-module algebras instead of
that of $\mathbb{C}[SL_n]_q$-comodule algebras mentioned above. This is what
we do in the present paper.

Let us recall now what the terminology "$U_q\mathfrak{sl}_n$-module algebra" means.
Let $A$ be a Hopf algebra. A unital algebra $F$ is said to be an $A$-module
algebra if $F$ is an $A$-module, the unit of $F$ is $A$-invariant (which
means $\xi(1)=\varepsilon(\xi)\cdot1$ for any $\xi\in A$), and, finally, the
multiplication $F\otimes F\rightarrow F$ intertwines the $A$-actions (we
recall that for any $A$-modules $V_1$, $V_2$ their tensor product is endowed
with an $A$-module structure via the comultiplication $\Delta:A\to A\otimes
A$).

\medskip

{\bf Remark.} In the sequel, we shall sometimes consider Hopf algebras with
an additional structure, namely, Hopf $*$-algebras (a Hopf $*$-algebra is a
pair $(A,*)$ where $A$ is a Hopf algebra and $*$ is an involution in $A$
with certain properties; see \cite{ChP}). In the case of Hopf $*$-algebras
the above-mentioned definition includes an additional requirement. Namely, let
$A_0=(A,*)$ be a Hopf $*$-algebra and $F$ an algebra. Then $F$ is said to be
an $A_0$-module algebra if, first, $F$ is an $A$-module algebra in the
previous sense, and, second, $F$ is involutive and the involutions in $A$
and $F$ agree as follows:
\begin{equation}\label{agree}
(\xi(f))^*=S(\xi)^*(f^*),\quad\xi\in A, f\in F.
\end{equation}

\medskip

(The notion of module algebras can be clarified in the
 classical setting
of a Lie group $G$ acting on a smooth $G$-space $X$. Denote by
$\mathfrak{g}$ the Lie algebra of $G$. Then the universal enveloping algebra
$U\mathfrak{g}$ acts on the space $C^\infty(X)$ via differential operators.
The usual Leibnitz rule means that $C^\infty(X)$ is a $U\mathfrak{g}$-module
algebra.)

Let us turn back to the quantum space of matrices. Now we are in position to
describe the very well known "infinitesimal version" of the left and right
actions of the quantum group $SL_n$ in $\mathbb{C}[M_{n}]_q$. Note, however,
that the left action we present below is not an analog of the classical one,
mentioned earlier. It is more convenient for us to use an action that
differs from the usual left one by a simple automorphism of
$U_q\mathfrak{sl}_n$.

\medskip

\begin{proposition}\label{uqslnxuqslnaction}

\noindent i) There exists a unique structure of $U_q\mathfrak{sl}_n$-module
algebra in $\mathbb{C}[M_{n}]_q$ such that
\begin{equation}\label{lh}
R(K_i)z_a^\alpha=\begin{cases}qz_a^\alpha, & a=i\\ q^{-1}z_a^\alpha, & a=i+1
\\ z_a^\alpha, & \mathrm{otherwise}\end{cases},
\end{equation}
\begin{equation}\label{lfe}
R(F_i)z_a^\alpha=\begin{cases}q^{1/2}z_{a+1}^\alpha, & a=i
\\ 0, & {\rm otherwise}\end{cases},\quad
R(E_i)z_a^\alpha=\begin{cases}q^{-1/2}z_{a-1}^\alpha, & a=i+1
\\ 0, &\mathrm{otherwise}\end{cases}.
\end{equation}
\noindent ii) There exists a unique structure of $U_q\mathfrak{sl}_n$-module
algebra in $\mathbb{C}[M_{n}]_q$ such that
\begin{equation}\label{rh}
L(K_j)z_a^\alpha=\begin{cases}qz_a^\alpha, & \alpha=n-j
\\ q^{-1}z_a^\alpha, &\alpha=n-j+1 \\ z_a^\alpha, &\mathrm{otherwise}
\end{cases},
\end{equation}
\begin{equation}\label{rfe}
L(F_j)z_a^\alpha=\begin{cases}q^{1/2}z_a^{\alpha+1}, & \alpha=n-j
\\ 0, &{\rm otherwise}\end{cases},\quad
L(E_j)z_a^\alpha=\begin{cases} q^{-1/2}z_a^{\alpha-1}, & \alpha=n-j+1
\\ 0, &\mathrm{otherwise}\end{cases}.
\end{equation}
\noindent iii) For any $\xi,\eta\in U_q\mathfrak{sl}_n$ the endomorphisms
$R(\xi)$, $L(\eta)$ commute $$R(\xi)L(\eta)f=L(\eta)R(\xi)f, \quad
f\in\mathbb{C}[M_{n}]_q.$$
\end{proposition}

\medskip
Note that by statement iii) in the above proposition, the algebra
$\mathbb{C}[M_{n}]_q$ is acted upon by the tensor product
$U_q\mathfrak{sl}_n\otimes U_q\mathfrak{sl}_n$:
$$\xi\otimes\eta(f)=R(\xi)L(\eta)f.$$

One can check that the $q$-determinant $\mathrm{det}_q (\mathbf{z})$
(\ref{sln}) is invariant with respect to both left and right
$U_q\mathfrak{sl}_n$-actions, i.e.
$$R(\xi)\mathrm{det}_q(\mathbf{z})=L(\xi)\mathrm{det}_q(\mathbf{z})=
\varepsilon(\xi)\cdot\mathrm{det}_q(\mathbf{z})$$ for any $\xi\in
U_q\mathfrak{sl}_n$. Thus the formulas from Proposition
\ref{uqslnxuqslnaction} define left and right $U_q\mathfrak{sl}_n$-actions
in $\mathbb{C}[SL_{n}]_q$. By analogy with the classical case, one has the following proposition (see
\cite{ChP}):

\medskip

\begin{proposition}\label{lrdecomp}
The $U_q\mathfrak{sl}_n\otimes U_q\mathfrak{sl}_n$-module
$\mathbb{C}[SL_n]_q$ splits into direct sum of simple pairwise
non-isomorphic submodules whose lowest vectors are given via $q$-minors as
follows
$$(z_n^n)^{a_1}\left(\mathbf{z}_{\,\,\,\,\,\,\,\{n-1,n\}}^{\wedge 2
\{n-1,n\}}\right)^{a_2}\left(
\mathbf{z}_{\,\,\,\,\,\,\,\{n-2,n-1,n\}}^{\wedge 3
\{n-2,n-1,n\}}\right)^{a_3}\cdots
\left(\mathbf{z}_{\;\;\;\;\;\;\,\,\,\,\,\,\,\,\,\,\{2,\dots,n\}}^{\wedge
(n-1) \{2,\ldots,n\}}\right)^{a_{n-1}}.$$
\end{proposition}

We recall that the $q$-minors
are defined by
\begin{equation}\label{minors}
(\mathbf{z}^{\wedge k})_{\{a_1,a_2,\ldots,a_k
\}}^{\{\alpha_1,\alpha_2,\ldots,\alpha_k
\}}\stackrel{\mathrm{def}}{=}\sum_{s \in
S_k}(-q)^{l(s)}z_{a_1}^{\alpha_{s(1)}}z_{a_2}^{\alpha_{s(2)}}\cdots
z_{a_k}^{\alpha_{s(k)}}
\end{equation}
with $\alpha_1<\alpha_2<\ldots<\alpha_k$, $a_1<a_2<\ldots<a_k$, and $l(s)$
being the length of $s\in S_k$. In particular,
$\mathrm{det}_q(\mathbf{z})=(\mathbf{z}^{\wedge n})_{\{1,2,\ldots,n
\}}^{\{1,2,\ldots,n \}}.$

\medskip

Let us denote the tensor product $U_q\mathfrak{sl}_n\otimes
U_q\mathfrak{sl}_n$ with the canonical Hopf algebra structure by
$U_q(\mathfrak{sl}_n\times\mathfrak{sl}_n)$. Thus, $\mathbb{C}[M_n]_q$ is a
$U_q(\mathfrak{sl}_n\times\mathfrak{sl}_n)$-module algebra. It follows from
the definition of the quantum universal enveloping algebra
$U_q\mathfrak{sl}_n$ that there is an embedding of Hopf algebras
$U_q(\mathfrak{sl}_n\times\mathfrak{sl}_n)\hookrightarrow
U_q\mathfrak{sl}_{2n}$ determined by $$1\otimes E_i\mapsto E_i, \quad
1\otimes F_i\mapsto F_i, \quad 1\otimes K^{\pm1}_i\mapsto K^{\pm1}_i, \quad
i=1,\ldots n-1,$$
$$E_i\otimes1\mapsto E_{n+i}, \quad F_i\otimes1\mapsto F_{n+i}, \quad
K^{\pm1}_i\otimes1\mapsto K^{\pm1}_{n+i}, \quad i=1,\ldots n-1.$$ This is a
$q$-analog of the embedding $SL_n\times SL_n\hookrightarrow SL_{2n}$ given,
in the matrix realization, by $$(A,B)\mapsto\left(\begin{array}{cl} A &
\quad 0\\ 0 & \quad B\end{array}\right).$$ In the next subsection we shall
extend the above $U_q(\mathfrak{sl}_n\times\mathfrak{sl}_n)$-module algebra
structure in $\mathbb{C}[M_{n}]_q$ to a structure of
$U_q\mathfrak{sl}_{2n}$-module algebra.

\medskip

\subsection{A structure of $U_q\mathfrak{sl}_{2n}$-module algebra on
$\mathbb{C}[M_{n}]_q$}

In this subsection we describe a "hidden" $U_q\mathfrak{sl}_{2n}$-module
algebra structure in $\mathbb{C}[M_{n}]_q$. It was discovered in
\cite{SV1}. Its classical counterpart comes from an embedding of the matrix
space $M_{n}$ into the Grassmannian $\mathrm{Gr}_{n}(\mathbb{C}^{2n})$ as the
affine cell $U\subset\mathrm{Gr}_{n}(\mathbb{C}^{2n})$ defined by the
inequality $t\neq0$ with $t$ being a distinguished Pl\"ucker coordinate. A
$q$-version of the embedding is described in \cite{SV1}, Proposition 0.7
(see also Proposition 5.4 from \cite {SSV1}).

Let us turn to the quantum case. The following statement was proved in
\cite[Section 2]{SSV1}.

\medskip

\begin{proposition}\label{uqsl2naction}
There exists a unique $U_q\mathfrak{sl}_{2n}$-module algebra structure in
$\mathbb{C}[M_{n}]_q$ given on the Hopf subalgebra
$U_q(\mathfrak{sl}_n\times\mathfrak{sl}_n)$ by the formulas from Proposition
\ref{uqslnxuqslnaction} and on the remaining generators $K_n^{\pm1}$, $F_n$,
$E_n$ by
\begin{equation}\label{h}
K_nz_a^\alpha=\begin{cases}q^2z_a^\alpha, & a=n \;\&\;\alpha=n
\\ qz_a^\alpha, & a=n \;\&\;\alpha \ne n \quad{\rm or}\quad a \ne n \;\&\;
\alpha=n \\ z_a^\alpha, & {\rm otherwise}\end{cases},
\end{equation}
\begin{equation}\label{fe}
F_nz_a^\alpha=q^{1/2}\begin{cases}1, & a=n \;\& \;\alpha=n \\ 0, &{\rm
otherwise}\end{cases}, \quad E_nz_a^\alpha=-q^{1/2}
\begin{cases}q^{-1}z_a^mz_n^\alpha, & a \ne n \;\&\;\alpha \ne n \\
(z_n^m)^2,
& a=n \;\&\;\alpha=n \\ z_n^mz_a^{\alpha}, & {\rm otherwise}\end{cases}.
\end{equation}
\end{proposition}

Let us point out some straightforward but essential properties of this
$U_q\mathfrak{sl}_{2n}$-action in $\mathbb{C}[M_{n}]_q$. Denote by
$U_q\mathfrak{s}(\mathfrak{gl}_{n}\times\mathfrak{gl}_{n})$ the Hopf
subalgebra in $U_q\mathfrak{sl}_{2n}$ derived from
$U_q(\mathfrak{sl}_{n}\times\mathfrak{sl}_{n})$ by adding the generators
$K_n^{\pm1}$. Clearly, elements of
$U_q\mathfrak{s}(\mathfrak{gl}_{n}\times\mathfrak{gl}_{n})$ preserve the
natural $\mathbb{Z}_+$-grading in $\mathbb{C}[M_{n}]_q$ given by powers of
monomials. It is also obvious that the generators $F_n$, $E_n$ act in
$\mathbb{C}[M_{n}]_q$ as endomorphisms of degrees $-1$ and $1$,
respectively. All this may be derived also from the following convenient
description of the $\mathbb{Z}_+$-grading in $\mathbb{C}[M_{n}]_q$:
\begin{equation}\label{deg}
\mathrm{deg} f=N \Leftrightarrow \hat{K}f=q^{2N}f
\end{equation}
where $\hat{K}$ is the element of the center of
$U_q\mathfrak{s}(\mathfrak{gl}_{n}\times\mathfrak{gl}_{n})$ given by
\begin{equation}\label{center}
\hat{K}=(K_n)^n\cdot\prod_{j=1}^{n-1}(K_jK_{2n-j})^j.
\end{equation}

For computational purposes, it is important to understand the structure of
$\mathbb{C}[M_{n}]_q$ as a
$U_q\mathfrak{s}(\mathfrak{gl}_{n}\times\mathfrak{gl}_{n})$-module in
greater details. The following statement is a straightforward consequence of
Proposition \ref{lrdecomp}.

\medskip

\begin{proposition}\label{huashmidt}
The $U_q\mathfrak{s}(\mathfrak{gl}_{n}\times\mathfrak{gl}_{n})$-module
$\mathbb{C}[M_n]_q$ splits into direct sum of simple pairwise non-isomorphic
submodules $\mathbb{C}[M_n]^{(k_1,k_2,\ldots,k_n)}_q$, $k_1\geq
k_2\geq\ldots\geq k_n\geq0$, whose lowest vectors are given by
$$(z_n^n)^{k_1-k_2}\left(\mathbf{z}_{\,\,\,\,\,\,\,\{n-1,n\}}^{\wedge 2
\{n-1,n\}}\right)^{k_2-k_3}\left(
\mathbf{z}_{\,\,\,\,\,\,\,\{n-2,n-1,n\}}^{\wedge 3
\{n-2,n-1,n\}}\right)^{k_3-k_4}\cdots (\det \nolimits_q \mathbf{z})^{k_n}.$$
\end{proposition}

\medskip

In what follows, the above $U_q\mathfrak{sl}_{2n}$-action in
$\mathbb{C}[M_{n}]_q$ will be sometimes called 'the initial' one, in
contrast to a twisted action described in the next subsection.

Let us present another view on the above $U_q\mathfrak{sl}_{2n}$-action in
$\mathbb{C}[M_{n}]_q$. The point is that the corresponding classical
$U\mathfrak{sl}_{2n}$-action is well known in the theory of bounded
symmetric domains (see, for instance, \cite{Arazy}). In framework of this
theory, it is constructed as follows. The vector space $M_{n}$ contains the
so-called matrix ball (the boundary symmetric domain of type $I_{n,n}$) $$
\mathscr{D}=\{\mathbf{z}\in M_{n}\,|\,\mathbf{z}\mathbf{z}^*<1\} $$ (with
$*$ being the hermitian conjugation and $1$ the unit matrix). It is known
that the real simple Lie group $SU_{n,n}$ acts on $\mathscr{D}$ via
biholomorphic automorphisms, and $S(U_n\times U_n)\subset SU_{n,n}$ is the
isotropy subgroup of the center $0\in\mathscr{D}$. Thus elements of the
universal enveloping algebra $U\mathfrak{su}_{n,n}$, and hence elements of
its complexification $U\mathfrak{sl}_{2n}$, act on the space of holomorphic
functions on $\mathscr{D}$ via differential operators. These differential
operators have polynomial coefficients and, thus, preserve
$\mathbb{C}[M_{n}]$. The resulting $U\mathfrak{sl}_{2n}$-action in
$\mathbb{C}[M_{n}]$ is what we call the initial one. In framework of this
approach, the result of Proposition \ref{huashmidt} is just a $q$-analog of
the famous Hua-Schmid decomposition \cite{Arazy} whereas the quantum
enveloping algebra
$U_q\mathfrak{s}(\mathfrak{gl}_{n}\times\mathfrak{gl}_{n})$ itself is an
analog of the universal enveloping algebra of the complexified Lie algebra
of the isotropy subgroup $S(U_n\times U_n)$.

\medskip

\subsection{A twisted $U_q\mathfrak{sl}_{2n}$-action on
$\mathbb{C}[M_{n}]_q$}

In this subsection we introduce a one-parameter family $\pi_\lambda$,
$\lambda\in\mathbb{R}$, of $U_q\mathfrak{sl}_{2n}$-actions in
$\mathbb{C}[M_{n}]_q$ such that the initial $U_q\mathfrak{sl}_{2n}$-action,
defined in the previous subsection, corresponds to $\lambda=0$. In the
classical case the corresponding twisted $U\mathfrak{sl}_{2n}$-action
$\pi_\lambda$ for $\lambda\in\mathbb{Z}$ can be produced by trivializing the
homogeneous line bundle $\mathcal{O}(-\lambda)$ on the Grassmannian
$\mathrm{Gr}_{n}(\mathbb{C}^{2n})$ over the affine cell $U$. Namely, we
identify the space of polynomials on $M_{n}$ with the space of sections
$\Gamma(U,\mathcal{O}(-\lambda))$ by $f(\mathbf{z})\sim f(\mathbf{z})\cdot
t^{-\lambda}$ (here $t$ is the distinguished Pl\"ucker coordinate mentioned at
the beginning of the previous subsection) and define the
$U\mathfrak{sl}_{2n}$-action $\pi_\lambda$ as follows:
\begin{equation}\label{altdes}
(\pi_\lambda(\xi)f)\cdot t^{-\lambda}=\xi(f\cdot t^{-\lambda}),\quad\xi\in
U\mathfrak{sl}_{2n}.
\end{equation}
Note that, among the actions $\pi_\lambda$, the initial action $\pi_0$ is
the only one that makes $\mathbb{C}[M_{n}]$ into a
$U\mathfrak{sl}_{2n}$-{\it module algebra}. This is true in the $q$-setting
as well.

Let us define a quantum version of the $U\mathfrak{sl}_{2n}$-action
(\ref{altdes}).

\medskip

\begin{proposition}
For any $\lambda\in\mathbb{R}$ the formulas $$\pi_ \lambda(K_j^{\pm 1})f=
\begin{cases} K_j^{\pm 1}f, & j \ne n \\ q^{\pm \lambda}K_n^{\pm 1}f, & j=n
\end{cases},$$
$$\pi_ \lambda(F_j)f=
\begin{cases} F_jf, & j \ne n \\ q^{-\lambda} F_nf, & j=n
\end{cases},\quad \pi_ \lambda(E_j)f=
\begin{cases} E_jf, & j \ne n \\ E_nf-q^{1/2}\dfrac{1-q^{2
\lambda}}{1-q^2}(K_nf)z_n^n, & j=n
\end{cases}$$
define a $U_q\mathfrak{sl}_{2n}$-action in $\mathbb{C}[M_{n}]_q$ (in the
right-hand sides the initial $U_q\mathfrak{sl}_{2n}$-action is used).
\end{proposition}
This proposition was proved in \cite[Proposition 6.2]{SSV2}. Note that
$\pi_0$ coincides with the initial $U_q\mathfrak{sl}_{2n}$-action. For
brevity, we shall denote the $U_q\mathfrak{sl}_{2n}$-module, corresponding
to $\lambda$, by $\mathcal{P}_\lambda$, namely
$\mathcal{P}_\lambda=(\mathbb{C}[M_{n}]_q,
U_q\mathfrak{sl}_{2n}, \pi_\lambda)$.

The classical counterpart of the above twisted
$U_q\mathfrak{sl}_{2n}$-action is also well known in the theory of bounded
symmetric domains. The corresponding $SU_{n,n}$-action (more precisely, the
action of the universal covering $\widetilde{SU}_{n,n}$) is defined by
\begin{equation}\label{weightaction}
\pi_\lambda(g):f(\mathbf{z})\mapsto f(g^{-1}\mathbf{z})\cdot
J_{g^{-1}}(\mathbf{z})^{\frac{\lambda}{2n}}
\end{equation}
with $J_{g^{-1}}(\mathbf{z})$ being the Jacobian of the biholomorphic map
$\mathbf{z}\mapsto g^{-1}\mathbf{z}$ (see \cite{Arazy} for details). For
$\lambda >2n-1$ the action $\pi_{\lambda}$ on a weighted Bergman space defines
a holomorphic discrete series representation of $\widetilde{SU}_{n,n}$. In
the last section of the paper we will describe a unitary structure on
$\mathcal{P}_\lambda$ which formally tends to the classical setting as $q\to
1$.

\bigskip

\section{Some $q$-differential operators}

\medskip

\subsection{Basic definitions}

One of our results is connected with a $q$-analog of "the wave operator"
$$\Box=\det\left(\frac{\partial}{\partial z_a^\alpha}\right).$$
We start with some general consideration
of $q$-differential
operators with constant coefficients.

To produce $q$-analogs of the partial derivatives, we use certain known
first order differential calculus over $\mathbb{C}[M_{n}]_q$, see \cite{ChP}. Let
$\Omega^1(M_{n})_q$ be the $\mathbb{C}[M_{n}]_q$-bimodule given by its
generators $dz_a^\alpha$, $a,\alpha=1,\ldots n$, and the relations $$
z_b^\beta dz_a^\alpha=\sum_{\alpha',\beta'=1}^n \sum_{a',b'=1}^n
R_{\beta\alpha}^{\beta'\alpha'}R^{b'a'}_{ba}dz_{a'}^{\alpha'}\cdot
z_{b'}^{\beta'},$$ with $$R^{b'a'}_{ba}=\begin{cases}q^{-1}, & a=b=a'=b' \\
1, & a \ne b \quad \&\quad a=a'\quad \&\quad b=b'\\ q^{-1}-q, & a<b \quad
\&\quad a=b' \quad \&\quad b=a'\\ 0, & {\rm otherwise}\end{cases}.$$ The map
$d:z_a^\alpha\mapsto dz_a^\alpha$ can be extended to a linear operator
$d:\mathbb{C}[M_{n}]_q\rightarrow\Omega^1(M_{n})_q$ satisfying the Leibnitz
rule $d(f_1f_2)=d(f_1)f_2+f_1d(f_2)$. The pair $(\Omega^1(M_{n})_q, d)$ is
the first order differential calculus over $\mathbb{C}[M_{n}]_q$ we need.

The calculus itself has been known for a long time \cite{ChP}. However, its hidden
$U_q\mathfrak{sl}_{2n}$-symmetry, was observed much later in \cite{SV1}. To
be more precise, it is proved in \cite{SV1} that there exists a unique
structure of a $U_q\mathfrak{sl}_{2n}$-module in $\Omega^1(M_{n})_q$ such
that, first, the map $d$ is a morphism of $U_q\mathfrak{sl}_{2n}$-modules,
and, second, the left and right multiplications
$$\mathbb{C}[M_{n}]_q\otimes\Omega^1(M_{n})_q\to\Omega^1(M_{n})_q,
\quad \Omega^1(M_{n})_q\otimes\mathbb{C}[M_{n}]_q\to\Omega^1(M_{n})_q$$ are
morphisms of $U_q\mathfrak{sl}_{2n}$-modules. This is usually expressed by
saying that the first order differential calculus $(\Omega^1(M_{n})_q, d)$
is $U_q\mathfrak{sl}_{2n}$-covariant. Before \cite{SV1} appeared, only
$U_q(\mathfrak{sl}_{n}\times\mathfrak{sl}_{n})$-covariance of the calculus
was known.

The first order differential calculus allows us to define the $q$-analogs of
partial derivatives as follows: Set $$d f=\sum_{a=1}^{n}\sum_{\alpha=1}^{n}
\frac{\partial f}{\partial z_a^{\alpha}}\cdot dz_a^{\alpha},\qquad f\in
\mathbb{C}[M_{n}]_q.$$ Here the left-hand side defines the right-hand one.

It is quite reasonable to regard the unital subalgebra in
$\mathrm{End}(\mathbb{C}[M_{n}]_q)$ generated by all the derivatives as an
analog of the algebra of differential operators with constant coefficients.
This algebra seems to be interesting in itself. First of all, it admits a very
explicit description. Namely, it is observed in \cite[Section 2]{Gauss} that
the map $z_a^\alpha\mapsto\frac{\partial}{\partial z^\alpha_a}$ may be
extended to an algebra homomorphism $\Upsilon:
\mathbb{C}[M_{n}]_q\to\mathrm{End}(\mathbb{C}[M_{n}]_q)$ which means that
the operators $\frac{\partial}{\partial z_a^{\alpha}}$ satisfy the same
commutation relations as the generators $z_a^{\alpha}$ of
$\mathbb{C}[M_n]_q$ do. Further, the algebra is invariant with respect to
a certain natural $U_q(\mathfrak{sl}_{n}\times\mathfrak{sl}_{n})$-action in
${\rm End}(\mathbb{C}[M_{n}]_q)$ defined via a $q$-analog of the commutator.
Let us describe this latter observation in full details.

Endow the space ${\rm End}(\mathbb{C}[M_{n}]_q)$ with a structure of
$U_q(\mathfrak{sl}_{n}\times\mathfrak{sl}_{n})$-module as follows: For
$\xi\in U_q(\mathfrak{sl}_{n}\times\mathfrak{sl}_{n})$, $T\in{\rm
End}(\mathbb{C}[M_{n}]_q)$ put $$\xi(T)=\sum_j\xi''_j\cdot T\cdot
S^{-1}(\xi'_j),$$ where $\sum_j\xi'_j\otimes\xi''_j=\Delta(\xi)$ (here
$\Delta$ denotes the comultiplication in
$U_q(\mathfrak{sl}_{n}\times\mathfrak{sl}_{n})$), $S$ is the antipode of
$U_q(\mathfrak{sl}_{n}\times\mathfrak{sl}_{n})$, and the elements in the
right-hand side are multiplied within ${\rm End}(\mathbb{C}[M_{n}]_q)$. It
is explained in \cite{Gauss} that the
$U_q(\mathfrak{sl}_{n}\times\mathfrak{sl}_{n})$-covariance of the first
order differential calculus $\left(\Omega^1(M_{n})_q, d\right)$ and the
explicit formulas for the
$U_q(\mathfrak{sl}_{n}\times\mathfrak{sl}_{n})$-action in
$\mathbb{C}[M_n]_q$, presented in Proposition \ref{uqslnxuqslnaction}, allow
one to prove $U_q(\mathfrak{sl}_{n}\times\mathfrak{sl}_{n})$-invariance of
the linear span of all $\frac{\partial}{\partial z^\alpha_a}$ in ${\rm
End}(\mathbb{C}[M_{n}]_q)$ and to describe the
$U_q(\mathfrak{sl}_{n}\times\mathfrak{sl}_{n})$-action on the partial
derivatives explicitly. The explicit description is based on the following
intertwining property of the homomorphism $\Upsilon$
$$ \Upsilon(\xi f)=\omega(\xi)\Upsilon(f),\quad \forall \xi\in
U_q(\mathfrak{sl}_{n}\times\mathfrak{sl}_{n}), \forall
f\in\mathbb{C}[M_{n}]_q$$ with $\omega$ being the automorphism of
$U_q(\mathfrak{sl}_{n}\times\mathfrak{sl}_{n})$ (the Chevalley involution)
given by $$ \omega(E_i)=-F_i,\quad \omega(F_i)=-E_i, \quad
\omega(K^{\pm1}_i)=K^{\mp1}_i.$$

\medskip

\subsection{A $q$-wave operator}

Here we apply the results from the previous subsection to study the
$q$-wave operator given by
$$
\Box_q=\sum_{s \in S_n}(-q)^{l(s)}\cdot\frac{\partial}{\partial
z_1^{s(1)}}\cdot\frac{\partial}{\partial z_2^{s(2)}}\cdot \ldots \cdot
\frac{\partial}{\partial z_n^{s(n)}}.$$ Clearly, the $q$-wave operator
belongs to the center of the algebra of quantum differential operators with
constant coefficients since $\Box_q=\Upsilon(\det \nolimits_q
(\mathbf{z}))$. The latter formula, together with
$U_q(\mathfrak{sl}_{n}\times\mathfrak{sl}_{n})$-invariance of the
$q$-determinant and the above intertwining property of $\Upsilon$, implies
also that the operator $\Box_q$ commutes with the action of
$U_q(\mathfrak{sl}_{n}\times\mathfrak{sl}_{n})$ in $\mathbb{C}[M_{n}]_q$.
Also, we can easily prove that
\begin{equation}\label{knbox}
K_n\cdot\Box_q=q^{-2}\Box_q\cdot K_n.
\end{equation}
Indeed, the degree of the operator $\Box_q$ in $\mathbb{C}[M_{n}]_q$ is
equal to $-n$ which means $\hat{K}\cdot\Box_q=q^{-2n}\cdot\Box_q\cdot
\hat{K}$ (see (\ref{deg})). The latter equality implies (\ref{knbox}) since
$\Box_q$ commutes with all the $K_i^{\pm1}$'s for $i\neq n$.

\bigskip

\section{A covariance property}

\medskip

\subsection{Formulation}

The intertwining properties of the $q$-wave operator derived above may be written
in a unified way as follows
\begin{equation}\label{obv}
\Box_q^{\,l}\cdot\pi_{n-l}(\xi)=\pi_{n+l}(\xi)\cdot\Box_q^{\,l},\qquad\xi\in
U_q\mathfrak{s}(\mathfrak{gl}_{n}\times\mathfrak{gl}_{n}).
\end{equation}
It turns out that this obvious symmetry of the operator $\Box_q^{\,l}$ is
a part of a large hidden symmetry: 

\medskip

\begin{theorem}\label{Bol}\hfill

\noindent For any $l\in\mathbb{N}$ the linear operator
$\Box_q^{\,l}:\mathbb{C}[M_{n}]_q\rightarrow \mathbb{C}[M_{n}]_q$
intertwines the $U_q\mathfrak{sl}_{2n}$-actions $\pi_{n-l}$ and $\pi_{n+l}$:
$$\Box_q^{\,l}\cdot\pi_{n-l}(\xi)=\pi_{n+l}(\xi)\cdot\Box_q^{\,l},\qquad\xi
\in
U_q\mathfrak{sl}_{2n}$$ (in other words, the map
$\Box_q^{\,l}:\mathcal{P}_{n-l}\to\mathcal{P}_{n+l}$ is a morphism of
$U_q\mathfrak{sl}_{2n}$-modules).
\end{theorem}

\medskip

We will proof Theorem \ref{Bol} in subsection 4.5. The proof uses some results from the theory of quantum bounded
symmetric domains which we recall in the subsequent three subsections. Very briefly, the idea is as follows (compare with
\cite{Arazy}): We use the $q$-Cauchy-Szeg\"o integral formula to rewrite the
operator $\Box_q^l$ as a $q$-integral operator; then, using some standard technique, we prove that the
$q$-integral operator intertwines the $U_q\mathfrak{sl}_{2n}$-actions
$\pi_{n-l}$ and $\pi_{n+l}$.

\medskip

\subsection {A $q$-analog of the Cauchy-Szeg\"o integral
representation}\label{sect}

For any bounded symmetric domain there is a multivariable generalization of
the famous Cauchy formula, the so-called Cauchy-Szeg\"o integral
representation \cite{Arazy, Hua}. This integral formula restores a
holomorphic function on the domain from its boundary value on the
Shilov boundary. In the case of the unit matrix ball the
Cauchy-Szeg\"o formula looks as follows

$$f(\mathbf{z})=\int
\limits_{S(\mathscr{D})}\frac{f(\boldsymbol{\zeta})}{\mathrm{det}
(1-\mathbf{z}\boldsymbol{\zeta}^*)^n}d \nu(\boldsymbol{\zeta}).$$ Here
$S(\mathscr{D})$ is the Shilov boundary of the unit matrix ball
$\mathscr{D}\in M_n$ $$ S(\mathscr{D})=\{\mathbf{z}\in M_n\,|\,
\mathbf{z}\mathbf{z}^*=1\},$$ and $d\nu$ is the unique ${U}_n$-invariant
normalized measure on $S(\mathscr{D})$ which, of course, coincides with the
Haar measure under the identification $S(\mathscr{D})=U_n$. 

A $q$-analog of this formula was found in \cite{CS} in framework of quantum
bounded symmetric domain theory. Particularly, in that paper $q$-analogs of
the Shilov boundary $S(\mathscr{D})$, the measure $d\nu$, and the kernel
$\mathrm{det}(1-\mathbf{z}\boldsymbol{\zeta}^*)^{-n}$ were found. In this
subsection we recall all these results. We omit proofs. An interested reader might 
want to look into \cite{CS} which is the main reference for this section.

The $q$-analog of the Shilov boundary is described by a (noncommutative) $*-algebra$ of functions on it. 
It is also natural to require the quantum Shilov boundary to be a homogeneous space of 
the $quantum$ group $SU_{n,n}$. Here is an explicit construction.

The localization of $\mathbb{C}[M_n]_q$ with respect to the multiplicative
system $\mathrm{det}_q (\mathbf{z})^\mathbb{N}$ is called the algebra of
regular functions on the quantum $GL_n$ and is denoted by
$\mathbb{C}[GL_n]_q$ (see subsection \ref{qsom}). It was observed in
\cite[Lemma 2.1]{CS} that there exists a unique involution $*$ in
$\mathbb{C}[GL_n]_q$ such that
\begin{equation}\label{involshilov}(z_a^\alpha)^*=(-q)^{a+\alpha-2n}\mathrm{det}_q(\mathbf{z})^{-1}\cdot
\mathbf{z}^{\wedge(n-1)J_\alpha}_{\,\,\,\,\,\,\,\,\,\,\,\,\,\,\,\,\,\,\,\,
J_a},
\end{equation}
with $J_c\stackrel{\mathrm{def}}=\{1,2,\ldots,n \}\backslash\{c\}$ (here we
use the notation (\ref{minors})). The $*$-algebra
$\mathrm{Pol}(S(\mathscr{D}))_q=(\mathbb{C}[GL_n]_q,*)$ is a $q$-analog of
the polynomial algebra on the Shilov boundary of the matrix ball
$\mathscr{D}$. Note that
\begin{equation}\label{detdet}
\mathrm{det}_q (\mathbf{z})\mathrm{det}_q (\mathbf{z})^*=\mathrm{det}_q
(\mathbf{z})^*\mathrm{det}_q (\mathbf{z})=q^{-n(n-1)}.
\end{equation}

Let's describe a structure of homogeneous space of 
the quantum group $SU_{n,n}$ on the quantum Shilov boundary. 
Recall that the $q$-determinant $\det \nolimits_q
(\mathbf{z})$ belongs to the center of $\mathbb{C}[M_n]_q$ and is a
"relative invariant" with respect to the
$U_q\mathfrak{s}(\mathfrak{gl}_{n}\times\mathfrak{gl}_{n})$-action:
\begin{equation}\label{relinv}
\xi\det \nolimits_q(\mathbf{z})=\varepsilon(\xi)\cdot\det \nolimits_q
(\mathbf{z}), \quad\xi\in
U_q(\mathfrak{sl}_{n}\times\mathfrak{sl}_{n}),\qquad K_n\det \nolimits_q
(\mathbf{z})=q^2\det \nolimits_q (\mathbf{z}).
\end{equation}
Using (\ref{relinv}), one can make $\mathrm{Pol}(S(\mathscr{D}))_q$ into a
$U_q\mathfrak{s}(\mathfrak{gl}_{n}\times\mathfrak{gl}_{n})$-module algebra. (More precisely, we can use the above formulas to define a $U_q\mathfrak{s}(\mathfrak{gl}_{n}\times\mathfrak{gl}_{n})$-action on negative powers of $\det \nolimits_q (\mathbf{z})$ which suffices to extend the $U_q\mathfrak{s}(\mathfrak{gl}_{n}\times\mathfrak{gl}_{n})$-action from $\mathbb{C}[M_n]_q$ to $\mathrm{Pol}(S(\mathscr{D}))_q$.)

In fact (see \cite[Section 2]{CS}), the above
$U_q\mathfrak{s}(\mathfrak{gl}_{n}\times\mathfrak{gl}_{n})$-module algebra
structure in $\mathrm{Pol}(S(\mathscr{D}))_q$ may be extended to a structure
of $U_q\mathfrak{sl}_{2n}$-module algebra which coincides on the subspace
$\mathbb{C}[M_n]_q\subset\mathrm{Pol}(S(\mathscr{D}))_q$ with the
$U_q\mathfrak{sl}_{2n}$-module algebra structure described in Proposition
\ref{uqsl2naction}.

Let us recall the definition of the "real form" $U_q\mathfrak{su}_{n,n}$ of
the quantum universal enveloping algebra $U_q\mathfrak{sl}_{2n}$.
$U_q\mathfrak{su}_{n,n}$ is simply the pair $(U_q\mathfrak{sl}_{2n},*)$ with
$*$ being an involution in $U_q\mathfrak{sl}_{2n}$ determined by
\begin{align*}
E_n^*&=-K_nF_n,&F_n^*&=-E_nK_n^{-1},&(K_n^{\pm 1})^*&=K_n^{\pm 1},&
\\ E_j^*&=K_jF_j,&F_j^*&=E_jK_j^{-1},&(K_j^{\pm 1})^*&=K_j^{\pm1},
&\text{for}\quad j \ne n.
\end{align*}
It is not difficult to verify that
$U_q\mathfrak{su}_{n,n}=(U_q\mathfrak{sl}_{2n},*)$ is a Hopf $*$-algebra
(see \cite{ChP} for definitions). Evidently, the involution $*$ keeps the Hopf subalgebras
$U_q(\mathfrak{sl}_{n}\times\mathfrak{sl}_{n})$ and
$U_q\mathfrak{s}(\mathfrak{gl}_{n}\times\mathfrak{gl}_{n})$ invariant, and
we shall denote the corresponding Hopf $*$-subalgebras in
$U_q\mathfrak{su}_{n,n}$ by $U_q(\mathfrak{su}_{n}\times\mathfrak{su}_{n})$
and $U_q\mathfrak{s}(\mathfrak{u}_{n}\times\mathfrak{u}_{n})$, respectively.

The crucial property of the involution (\ref{involshilov}) is the following observation: it
makes $\mathrm{Pol}(S(\mathscr{D}))_q$ into a
$U_q\mathfrak{su}_{n,n}$-module algebra (this is explained in \cite{CS}
after Proposition 2.7). It is in this sense that the quantum Shilov boundary is a
homogeneous space of the quantum group $SU_{n,n}$.

To define a $q$-analog of the measure $d\nu$ on $S(\mathscr{D})$, we note
(see \cite[Section 3]{CS}) that the $*$-algebra
$\mathrm{Pol}(S(\mathscr{D}))_q$ is closely related to the $*$-algebra
$\mathbb{C}[{U}_n]_q=(\mathbb{C}[GL_n]_q,\star)$ of regular functions on the
quantum group ${U}_n$ where, we recall, the involution $\star$ is
defined by $(z_a^\alpha)^\star=(-q)^{a-\alpha}(\det \nolimits_q
\mathbf{z})^{-1}\cdot
\mathbf{z}^{\wedge(n-1)J_\alpha}_{\,\,\,\,\,\,\,\,\,\,\,\,\,\,\,\,\,\,\,\,J_a}
$. \footnote{The quantum group ${U}_n$ is one of the most well studied objects in quantum group theory. We refer to \cite{Koe} for basic definitions and facts about this quantum group. Of course, there are many other good references.}

It's easy to check that $*=\theta^{-1}\cdot\star\cdot \theta$ where
$\theta:\mathbb{C}[GL_n]_q \to \mathbb{C}[GL_n]_q$ is an automorphism given
by $\theta:z_a^\alpha \mapsto q^{\alpha-n}z_a^\alpha$. It is  known
that the compact quantum group $U_n$ possesses a unique
normalized invariant integral, an analog of the Haar integral. The
isomorphism $\theta$ of $*$-algebras $\mathrm{Pol}(S(\mathscr{D}))_q \to
\mathbb{C}[{U}_n]_q$ allows us to "transfer" the invariant integral onto
$\mathrm{Pol}(S(\mathscr{D}))_q$. In this way we get a positive
$U_q\mathfrak{s}(\mathfrak{gl}_{n}\times\mathfrak{gl}_{n})$-invariant linear
functional $\mathrm{Pol}(S(\mathscr{D}))_q \to \mathbb{C}$, $f\mapsto
\int_{S(\mathscr{D})_q}fd \nu$ which is the analog of the integral with
respect to $d\nu$. The
$U_q\mathfrak{s}(\mathfrak{gl}_{n}\times\mathfrak{gl}_{n})$-invariance means
\begin{equation}\label{inva}
\int\limits_{S(\mathscr{D})_q}\!\!\! \xi f\,d\nu=
\varepsilon(\xi)\cdot\int\limits_{S(\mathscr{D})_q}\!\!\! f\,d\nu,\quad
\forall\xi\in U_q\mathfrak{s}(\mathfrak{gl}_{n}\times\mathfrak{gl}_{n}).
\end{equation}

Finally let us describe the analog of the Cauchy-Szeg\"o kernel
$\mathrm{det}(1-\mathbf{z}\boldsymbol{\zeta}^*)^{-n}$.

Consider the algebra $\mathrm{Pol}(M_n \times
S(\mathscr{D}))_q=\mathbb{C}[M_n]_q^{\mathrm{op}}\otimes
\mathrm{Pol}(S(\mathscr{D}))_q$ with "op" indicating the change of the
multiplication to the opposite one. Equip it with a $\mathbb{Z}_+$-grading
by setting $ \deg(z_a^\alpha \otimes f)=1$ for any
$f\in\mathrm{Pol}(S(\mathscr{D}))_q$. Its completion with respect to this
grading is denoted by $\mathrm{Fun}(M_n \times S(\mathscr{D}))_q$. The
elements of $\mathrm{Fun}(M_n \times S(\mathscr{D}))_q$ are $q$-analogs of
kernels of integral operators, while the elements of the subalgebra
$\mathrm{Pol}(M_n \times S(\mathscr{D}))_q$ are $q$-analogs of polynomial
kernels.

Let us comment on the replacement of the multiplication law in the first
tensor multiplier in the definition of the algebra $\mathrm{Pol}(M_n \times
S(\mathscr{D}))_q$. Given a Hopf algebra $A$ and two $A$-module algebras
$F_1, F_2$, $A$-invariant elements in $F_1\otimes F_2$ {\it do not} form a
subalgebra. However, they do form a subalgebra in $F_1^{\mathrm{op}}\otimes
F_2$ \cite{V1}. Almost all the kernels we encounter in the present paper are
$U_q\mathfrak{s}(\mathfrak{gl}_{n}\times\mathfrak{gl}_{n})$-invariant in
$\mathrm{Pol}(M_n \times S(\mathscr{D}))_q$ or $\mathrm{Fun}(M_n \times
S(\mathscr{D}))_q$ and so, as we have explained, form a subalgebra.

Let us explain now why we are interested in
$U_q\mathfrak{s}(\mathfrak{gl}_{n}\times\mathfrak{gl}_{n})$-invariant
kernels. It is not difficult to prove that there is a one-to-one
correspondence between
$U_q\mathfrak{s}(\mathfrak{gl}_{n}\times\mathfrak{gl}_{n})$-invariant
elements in $\mathrm{Fun}(M_n \times S(\mathscr{D}))_q$ and endomorphisms of
the $U_q\mathfrak{s}(\mathfrak{gl}_{n}\times\mathfrak{gl}_{n})$-module
$\mathbb{C}[M_n]_q$, explicitly given as follows: the element
$K\in\mathrm{Fun}(M_n \times S(\mathscr{D}))_q$ defines the morphism
$f\mapsto(1\otimes\int_{S(\mathscr{D})_q})(K\cdot(1\otimes f))$. In other words, a
linear operator on $\mathbb{C}[M_n]_q$ intertwines the
$U_q\mathfrak{s}(\mathfrak{gl}_{n}\times\mathfrak{gl}_{n})$-action if and
only if it is a $q$-integral operator with an invariant kernel.

It is convenient to choose some generators of $\mathrm{Pol}(M_n \times
S(\mathscr{D}))_q$ and express all other invariant kernels from
$\mathrm{Pol}(M_n \times S(\mathscr{D}))_q$ or $\mathrm{Fun}(M_n \times
S(\mathscr{D}))_q$ as (finite or formal) series in those generators.

Consider the elements $\chi_k\in\mathrm{Pol}(M_n \times S(\mathscr{D}))_q$,
$k=1,\ldots, n$, given by
\begin{equation}\label{chi}
\chi_k=\sum{\mathbf{z}^{\wedge k}}_{J''}^{J'}\otimes
\left({\mathbf{z}^{\wedge k}}_{J''}^{J'}\right)^*
\end{equation}
where the sum is taken over the pairs of subsets $J',J''\subset
\{1,2,\ldots,n\}$ of cardinality $k$. It turns out that the elements
$\chi_k$ are pairwise commuting and
$U_q\mathfrak{s}(\mathfrak{gl}_{n}\times\mathfrak{gl}_{n})$-invariant
\cite[Section 10]{SSV2}.

\medskip

\begin{proposition}\label{gener}
The elements $\chi_1,\ldots,\chi_n$ generate the subalgebra of
$U_q\mathfrak{s}(\mathfrak{gl}_{n}\times\mathfrak{gl}_{n})$-invariant
kernels in $\mathrm{Pol}(M_n \times S(\mathscr{D}))_q$ (which is therefore a
commutative algebra).
\end{proposition}
\noindent{\bf Sketch of a proof.} Recall (see Proposition \ref{huashmidt})
that the $U_q\mathfrak{s}(\mathfrak{gl}_{n}\times\mathfrak{gl}_{n})$-module
$\mathbb{C}[M_n]_q$ splits into direct sum of simple pairwise non-isomorphic
submodules $\mathbb{C}[M_n]^{\mathbf{k}}_q$ with
$\mathbf{k}=(k_1,k_2,\ldots,k_n)$, $k_1\geq k_2\geq\ldots\geq k_n\geq0$.
Thus, any endomorphism of the
$U_q\mathfrak{s}(\mathfrak{gl}_{n}\times\mathfrak{gl}_{n})$-module
$\mathbb{C}[M_n]_q$ is a (in general infinite) series of the form
$\sum_{\mathbf{k}}c_{\,\mathbf{k}}\cdot\mathbf{P}^{\mathbf{k}}$ where
$\mathbf{P}^{\mathbf{k}}$ stands for the projection in $\mathbb{C}[M_n]_q$
onto $\mathbb{C}[M_n]^{\mathbf{k}}_q$ parallel to the sum of other
${U_q\mathfrak{s}(\mathfrak{gl}_{n}\times\mathfrak{gl}_{n})}$-submodules and
$c_{\,\mathbf{k}}$ are complex numbers. It is sufficient to show that each
projection $\mathbf{P}^{\mathbf{k}}$ is a $q$-integral operator whose
kernel is a function of $\chi_1,\ldots,\chi_n$. This may be done by
using well known orthogonality relations for the quantum group $U_n$ (see
\cite{Koe}) and the precise relation between the quantum Shilov boundary and
the quantum $U_n$ described previously.\hfill$\blacksquare$

\medskip

The projection $\mathbf{P}^{\mathbf{k}}$
can be written
as a $q$-integral operator
with a kernel
${P^{\mathbf{k}}\in\mathrm{Pol}(M_n\times S(\mathscr{D}))_q}$.
Namely, let $u^{\mathbf{k}}$ be a
polynomial such that $P^{\mathbf{k}}(\mathbf{z},\boldsymbol{\zeta}^*)=
u^{\mathbf{k}}(\chi_1, \chi_2,\ldots, \chi_n)$. Consider the isomorphism
\begin{equation}\label{isom}
\mathrm{Fun}(M_n \times
S(\mathscr{D}))_q^{U_q\mathfrak{s}(\mathfrak{gl}_{n}\times\mathfrak{gl}_{n})}
\to\mathbb{C}[[x_1,x_2,\ldots,x_n]]^{S_n},\quad \chi_k\mapsto\sigma_k, \quad
k=1,2,\ldots n
\end{equation}
from the subalgebra $\mathrm{Fun}(M_n \times
S(\mathscr{D}))_q^{U_q\mathfrak{s}(\mathfrak{gl}_{n}\times\mathfrak{gl}_{n})}$
of $U_q\mathfrak{s}(\mathfrak{gl}_{n}\times\mathfrak{gl}_{n})$-invariants in
$\mathrm{Fun}(M_n \times S(\mathscr{D}))_q$ to the algebra of symmetric
formal series of the variables $x_1,x_2,\ldots,x_n$, where $\sigma_k$ is the
$i$-th elementary symmetric polynomial in $x_1,x_2,\ldots,x_n$. The
image of $u^{\mathbf{k}}(\chi_1, \chi_2,\ldots, \chi_n)$ under the above
isomorphism differs only by a constant from the so-called Schur polynomial
$s_{\mathbf k}$ associated to the partition $\mathbf{k}$ (see \cite{M}),
viz.
\begin{equation}\label{COEF}
u^{\mathbf{k}}(\sigma_1, \sigma_2,\ldots, \sigma_n)=C(\mathbf{k})\cdot
s_{\mathbf{k}}(x_1,x_2,\ldots,x_n).
\end{equation}
We compute the coefficients $C(\mathbf{k})$ in the next subsection.

From now on, for a kernel
$K\in\mathrm{Fun}(M_n \times S(\mathscr{D}))_q$ we shall sometimes write
$$K(\mathbf{z},\boldsymbol{\zeta}^*),\quad
\int_{S(\mathscr{D})_q}K(\mathbf{z},\boldsymbol{\zeta}^*)\cdot
f(\boldsymbol{\zeta})d\nu(\boldsymbol{\zeta}),$$ instead of $K$ and
$(1\otimes\int_{S(\mathscr{D})_q})(K\cdot(1\otimes f))$, respectively.

Now we are ready to present the $q$-analog of the Cauchy-Szeg\"o integral
formula found in \cite{CS}. In short, it represents the identity operator on
$\mathbb{C}[M_n]_q$ in the form of a $q$-integral operator.
\begin{theorem}\label{CSK}\cite[Section 5]{CS} For any element $f \in
\mathbb{C}[M_n]_q$ one has $$ f(\mathbf{z})=\int
\limits_{S(\mathscr{D})_q}C_q(\mathbf{z},\boldsymbol{\zeta}^*)
f(\boldsymbol{\zeta})d \nu(\boldsymbol{\zeta})$$ where $C_q=
\prod\limits_{j=0}^{n-1}\left(1+\sum\limits_{k=1}^n(-q^{2j})^k
\chi_k\right)^{-1}$ (a $q$-analog of the Cauchy-Szeg\"o kernel).
\end{theorem}

\medskip

\subsection{$q$-Analogs of the kernels $\mathrm{det}
(1-\mathbf{z}\boldsymbol{\zeta}^*)^{-N}$}\label{kernels}

Consider the family $K_N$, $N=1,2,\ldots$, of
$U_q\mathfrak{s}(\mathfrak{gl}_{n}\times\mathfrak{gl}_{n})$-invariant
kernels given by
\begin{equation}\label{kern}
K_N=\prod_{j=0}^{N-1}\left(1+\sum_{k=1}^n(-q^{2j})^k \chi_k\right)^{-1}.
\end{equation}
Clearly, the $q$-Cauchy-Szeg\"o kernel $C_q$ defined in the previous
subsection coincides with $K_n$. One also has
$$\lim_{q\rightarrow1}K_N(\mathbf{z},\boldsymbol{\zeta}^*)=
\det(1-\mathbf{z}\boldsymbol{\zeta}^*)^{-N}$$ (the limit should be
understood formally). The aim of this subsection is to study these kernels
and the associated $q$-integral operators in details.

Let $\mathbf{K}_N$ be the
$U_q\mathfrak{s}(\mathfrak{gl}_{n}\times\mathfrak{gl}_{n})$-intertwining
$q$-integral operator corresponding to the kernel $K_N$:
$$\mathbf{K}_Nf(\mathbf{z})=\int_{S(\mathscr{D})_q}
K_N(\mathbf{z},\boldsymbol{\zeta}^*)\cdot
f(\boldsymbol{\zeta})d\nu(\boldsymbol{\zeta}).$$ Then
\begin{equation}\label{decom}
\mathbf{K}_N=\sum_{\mathbf{k}}c^{\mathbf{k}}_N\cdot \mathbf{P}^{\mathbf{k}}
\end{equation}
where $\mathbf{P}^{\mathbf{k}}$ stands for the projection defined in the
previous subsection (see the proof of Proposition \ref{gener}).
We are interested in an explicit formula for the coefficients
$c^{\mathbf{k}}_N$.

One may write (\ref{decom}) as an equality of kernels:
\begin{equation}\label{decomker}
K_N(\mathbf{z},\boldsymbol{\zeta}^*)=\sum_{\mathbf{k}}c^{\mathbf{k}}_N\cdot
P^{\mathbf{k}}(\mathbf{z},\boldsymbol{\zeta}^*).
\end{equation}
This approach, along with the formula (\ref{COEF}), allows us to use some
identities for the Schur functions \cite{M} to compute $c^{\mathbf{k}}_N$.

Our first step towards computing the coefficients in (\ref{decom}) consists
in computing the constants $C(\mathbf{k})$ in (\ref{COEF}). For that
purpose, we note that
$$K_n(\mathbf{z},\boldsymbol{\zeta}^*)=\sum_{\mathbf{k}}
P^{\mathbf{k}}(\mathbf{z},\boldsymbol{\zeta}^*)$$ (this is just another way
to formulate Theorem \ref{CSK}). In view of the explicit form of the
$q$-Cauchy-Szeg\"o kernel (Theorem \ref{CSK}) and the isomorphism
(\ref{isom}), the latter can be written as follows
$$\prod_{j=0}^{n-1}\left(1+\sum_{k=1}^n(-q^{2j})^k
\sigma_k\right)^{-1}=\sum_{\mathbf{k}} u^{\mathbf{k}}(\sigma_1,
\sigma_2,\ldots, \sigma_n),$$ or,
by taking into account (\ref{COEF})
\begin{equation}
\label{coef-a}
\prod_{j=0}^{n-1}\left(1+\sum_{k=1}^n(-q^{2j})^k
\sigma_k\right)^{-1}=\sum_{\mathbf{k}}C(\mathbf{k})\cdot
s_{\mathbf{k}}(x_1,x_2,\ldots,x_n).
\end{equation}
Recall, for any integer $N\ge 0$,
 the $q$-Pochhammer symbol
$(x;q^2)_N=\prod_{j=0}^{N-1}(1-xq^{2j})$.
We have then,
$$
\prod_{j=0}^{N-1}\left(1+\sum_{k=1}^n(-q^{2j})^k
\sigma_k\right)^{-1}
=\prod_{j=0}^{N-1}\left(\prod_{i=1}^n (1-q^{2j}x_i)
\right)^{-1}
=\prod_{i=1}^{n}\frac1{(x_i;q^2)_N};
$$
in particular for $N=n$ the
equality (\ref{coef-a}) reads
$$\prod_{i=1}^{n}\frac1{(x_i;q^2)_n}=\sum_{\mathbf{k}}C(\mathbf{k})\cdot
s_{\mathbf{k}}(x_1,x_2,\ldots,x_n).$$
Now we are in position to make use of
the following formula from \cite{M}
\begin{equation}\label{milne}
\prod_{i=1}^{n}\frac{(ax_i;q^2)_\infty}{(x_i;q^2)_\infty}=\sum_{\mathbf{k}}
C(\mathbf{k};a)\cdot
s_{\mathbf{k}}(x_1,x_2,\ldots,x_n)
\end{equation}
where
\begin{equation}\label{milne1}
C(\mathbf{k};a)=\prod_{i=1}^n\frac{(aq^{2-2i};q^2)_{k_i}\cdot
q^{2(i-1)k_i}}{(q^2;q^2)_{k_i+n-i}}\cdot\prod_{1\leq i<j\leq
n}(1-q^{2k_i-2k_j-2i+2j}).
\end{equation}
In our case $a=q^{2n}$ and thus the equality (\ref{COEF})
acquires the form
\begin{equation}\label{COEF1}
u^{\mathbf{k}}(\sigma_1, \sigma_2,\ldots,
\sigma_n)=C(\mathbf{k};q^{2n})\cdot s_{\mathbf{k}}(x_1,x_2,\ldots,x_n),
\end{equation}
where
$$C(\mathbf{k};q^{2n})
=\prod_{i=1}^n\frac{q^{2(i-1)k_i}}{(q^2;q^2)_{n-i}}\cdot\prod_{1\leq i<j\leq
n}(1-q^{2k_i-2k_j-2i+2j}).$$

We turn back now to computing the coefficients in (\ref{decom}). Identifying
$K_N$ with its image under the isomorphism (\ref{isom}), we have, in view of
(\ref{kern}),
$$
K_N=
\prod_{j=0}^{N-1}\left(1+\sum_{k=1}^n(-q^{2j})^k
\sigma_k\right)^{-1}=\prod_{i=1}^{n}\frac1{(x_i;q^2)_N}.
$$
By (\ref{milne})
and (\ref{COEF1})
\begin{equation}\label{coef-b}
K_N=\sum_{\mathbf{k}}C(\mathbf{k};q^{2N})\cdot
s_{\mathbf{k}}(x_1,x_2,\ldots,x_n)=\sum_{\mathbf{k}}
\frac{C(\mathbf{k};q^{2N})}
{C(\mathbf{k};q^{2n})}u^{\mathbf{k}}(\sigma_1,\sigma_2,\ldots, \sigma_n).\end{equation}
We have thus obtained

\medskip

\begin{proposition}\label{coeffi} The coefficients in (\ref{decom}) are given
by
$$c^{\mathbf{k}}_N=\frac{C(\mathbf{k};q^{2N})}
{C(\mathbf{k};q^{2n})}=\prod_{i=1}^n\frac{(q^{2N+2-2i};q^2)_{k_i}}
{(q^{2n+2-2i};q^2)_{k_i}}.$$
\end{proposition}

\medskip

\subsection{A $q$-analog of the Fock inner product}\label{FFIP} The
aim of this subsection is to describe some results on a $q$-analog of the
Fock inner product in $\mathbb{C}[M_n]$  obtained in
\cite{Gauss}. At the end of the subsection we shall prove a $q$-analog of
one known result by J.~Faraut and A.~Koranyi \cite{FK} which compares the
Fock inner product with the one in the Hilbert space of
square-integrable functions on the Shilov boundary of the matrix ball.

Recall that the Fock inner product in the space $\mathbb{C}[M_n]$ is
defined by
\begin{equation}\label{clfip}
(f_1\,,\,f_2)_F=\int\limits_{M_n}f_1(\mathbf{z})
\overline{f_2(\mathbf{z})}e^{-\mathrm{tr}(\mathbf{z}\mathbf{z}^*)}d\mathbf{z}
\end{equation}
with $d\mathbf{z}$ being the Lebesgue measure on $M_n$ normalized so that
$\left(1\,,\,1\right)_F=1$. The inner product possesses the following
remarkable property
\begin{equation}\label{1pro}
(\frac{\partial f_1}{\partial z_a^\alpha}\,,\,f_2)_F= (f_1\,,\,z_a^\alpha
f_2)_F \quad \forall a,\alpha.
\end{equation}
This property, along with $S(U_n\times U_n)$-invariance of the inner
product, is quite useful in explicit computations of various norms.

Below we present a $q$-analog of the Fock inner product. But first
we have to explain what is understood by an invariance of an
inner product in the $q$-setting.

Let $A_0=(A,*)$ be a Hopf $*$-algebra. An inner product $(\,,\,)$ on an
$A$-module $V$ is said to be $A_0$-invariant if for all $v_1,v_2\in V$ and
any $\xi\in A$
$$(\xi v_1,v_2)=(v_1,\xi^*v_2).$$

The following is one of the main results of \cite{Gauss}.

\medskip

\begin{proposition}\label{FFI}
There exists a (unique)
$U_q\mathfrak{s}(\mathfrak{u}_{n}\times\mathfrak{u}_{n})$-invariant inner
product $(\,\cdot\,,\,\cdot\,)_{F}$ in $\mathbb{C}[M_n]_q$ satisfying the
properties
$$\left(1\,,\,1\right)_{F}=1,$$
$$(\frac{\partial f_1}{\partial z_a^\alpha}\,,\,f_2)_{F}= (f_1\,,\, f_2\cdot
z_a^\alpha)_{F} \quad \forall a,\alpha.$$
\end{proposition}

\medskip

Let $(\,\cdot\,,\,\cdot\,)_{S(\mathscr{D})}$ be the inner product in
$\mathbb{C}[M_n]_q$ defined via the
$U_q\mathfrak{s}(\mathfrak{gl}_{n}\times\mathfrak{gl}_{n})$-invariant
integral on the quantum Shilov boundary
$$(f_1,f_2)_{S(\mathscr{D})}=\int
\limits_{S(\mathscr{D})_q}f_2(\boldsymbol{\zeta})^*f_1(\boldsymbol{\zeta}) d
\nu(\boldsymbol{\zeta}).$$ Clearly, the inner product is
$U_q\mathfrak{s}(\mathfrak{u}_{n}\times\mathfrak{u}_{n})$-invariant. This is
a consequence of (\ref{inva}) and the condition (\ref{agree}). Since
$(\,\cdot\,,\,\cdot\,)_{F}$ and $(\,\cdot\,,\,\cdot\,)_{S(\mathscr{D})}$ are
$U_q\mathfrak{s}(\mathfrak{u}_{n}\times\mathfrak{u}_{n})$-invariant, the
subspaces $\mathbb{C}[M_n]^{\mathbf{k}}_q$ are pairwise orthogonal with
respect to both inner products, and the corresponding norms are proportional
by the Schur lemma. The proportionality constant
 is computed in the classical
setting by J.~Faraut and A.~Koranyi \cite[Corollary 3.5]{FK}. Here we
present a $q$-analog of their result.

\medskip

\begin{proposition}\label{qfk}
\begin{equation}\label{FF}
(f_1,f_2)_{F}=\frac{\prod_{i=1}^n(q^{2n+2-2i};q^2)_{k_i}}
{(1-q^2)^{k_1+k_2+\ldots+k_n}}\cdot(f_1,f_2)_{S(\mathscr{D})}, \quad
f_1,f_2\in\mathbb{C}[M_n]^{\mathbf{k}}_q.
\end{equation}
\end{proposition}
\noindent{\bf Proof.} To prove the proposition, we need an explicit
description of the inner product $(\,\cdot\,,\,\cdot\,)_{F}$. 

Consider the algebra $\mathbb{C}[M_n \times
M_n]_q=\mathbb{C}[M_n]_q\otimes\mathbb{C}[M_n]_q$. Equip it with the natural
bigrading by setting $ \deg(f_1\otimes f_2)=(\deg(f_1),\deg(f_2))$ for any
$f_1, f_2 \in\mathbb{C}[M_n]_q$. Its completion with respect to this
bigrading is denoted by $\mathbb{C}[[M_n \times M_n]]_q$.

Let
$$ \hat{\chi}_k=\sum_{\genfrac{}{}{0mm}{}{J',J''
\subset \{1,2,\ldots,n\}}
{\mathrm{card}(J')=\mathrm{card}(J'')=k}}{\mathbf{z}^{\wedge
k}}_{J''}^{J'}\otimes {\mathbf{z}^{\wedge k}}_{J''}^{J'}\in \mathbb{C}[M_n
\times M_n]_q,\quad k=1,\ldots, n. $$
Note that  $\hat{\chi}_k$ are
similar to the kernels $\chi_k$ defined in (\ref{chi}). Since the latter
pairwise commute, the elements $\hat{\chi}_k$ pairwise commute as well. Put
$$\hat{K}_{\infty}=\prod_{j=0}^{\infty}\left(1+\sum_{k=1}^n(-q^{2j})^k
\hat{\chi}_k\right)^{-1}\in\mathbb{C}[[M_n \times M_n]]_q. $$ Let
$\langle\,\cdot\,,\,\cdot\,\rangle$ be the inner product in
$\mathbb{C}[M_n]_q$ so that $\hat{K}_{\infty}$
is the reproducing kernel, namely,
writing $\hat{K}_{\infty}=\sum_jk'_j\otimes k''_j$ we have then
$$f=\sum_jk'_j\cdot\langle f\,,\,k''_j\rangle$$.

\medskip

\begin{lemma}
$$\langle\frac{\partial f_1}{\partial
z_a^\alpha}\,,\,f_2\rangle=\frac1{1-q^2}\cdot\langle f_1\,,\, f_2\cdot
z_a^\alpha\rangle \quad \forall a,\alpha.$$
\end{lemma}

\medskip

\noindent{\bf Sketch of a proof.} The inner product
$\langle\,\cdot\,,\,\cdot\,\rangle$ is described in a slightly different way
in \cite[Theorem 6.1]{Gauss}. The equivalence of the two definitions may be
deduced from Theorem 9.1 in \cite{SSV2} via the limit
$\lambda\to\infty$.\hfill$\blacksquare$.
\medskip

The above lemma allows us to express the inner product
$\langle\,\cdot\,,\,\cdot\,\rangle$ via the $q$-Fock one
\begin{equation}\label{relation}
\langle f_1\,,\, f_2\rangle=(1-q^2)^{k}\cdot(f_1\,,\, f_2)_{F}
\end{equation}
for $f_1,f_2$ homogeneous of degree $k$. Thus, to prove the theorem it
suffices to show that
$$
\langle f_1\,,\, f_2\rangle=\prod_{i=1}^n(q^{2n+2-2i};q^2)_{k_i}
\cdot(f_1,f_2)_{S(\mathscr{D})}, \quad
f_1,f_2\in\mathbb{C}[M_n]^{\mathbf{k}}_q.
$$
It follows from the $q$-Cauchy-Szeg\"o formula (Theorem \ref{CSK}) that the
inner product $(\,\cdot\,,\,\cdot\,)_{S(\mathscr{D})}$ is the one associated
to the kernel
$$\hat{K}_{n}=\prod_{j=0}^{n-1}\left(1+\sum_{k=1}^n(-q^{2j})^k
\hat{\chi}_k\right)^{-1}\in\mathbb{C}[[M_n \times M_n]]_q $$
in the same sense as described  above for
$\hat{K}_{\infty}$.
 It remains to
use the same arguments as in the previous subsection and to compare the
coefficients $C(\mathbf{k};0)$ and $C(\mathbf{k};q^{2n})$ (see
(\ref{milne1})).\hfill$\blacksquare$

\medskip

In the last section of the present paper we shall present more general
result which is  due to \O{}rsted
\cite{Orsted-tams-80} and J.~Faraut and A.~Koranyi \cite{FK}
in the classical
setting.

\medskip

\subsection{Proof of the covariance property}

Now we are in position to prove Theorem \ref{Bol}.

In view of Theorem \ref{CSK}, $\Box_q^{\,l}$ is the integral operator with
the kernel $\Box_q^{\,l}K_n(\mathbf{z},\boldsymbol{\zeta}^*)$ (here and
further we assume that the operator $\Box_q$ acts on a kernel in the first
argument $\mathbf{z}$). We are going to compute the
kernel explicitly.

Recall the notation $P^{\mathbf{k}}(\mathbf{z},\boldsymbol{\zeta}^*)$ for
the kernel of the $q$-integral operator $\mathbf{P}^{\mathbf{k}}:
\mathbb{C}[M_n]_q\to\mathbb{C}[M_n]^{\mathbf{k}}_q$ (see subsection
\ref{kernels})
$$\mathbf{P}^{\mathbf{k}}f(\mathbf{z})=\int
\limits_{S(\mathscr{D})_q}P^{\mathbf{k}}(\mathbf{z},\boldsymbol{\zeta}^*)
f(\boldsymbol{\zeta})d \nu(\boldsymbol{\zeta}). $$
By Theorem \ref{CSK},
$$C_q(\mathbf{z},\boldsymbol{\zeta}^*)=\sum_{\mathbf{k}}
P^{\mathbf{k}}(\mathbf{z},\boldsymbol{\zeta}^*).$$
 Thus to compute the
kernel $\Box_q^{\,l}C_q(\mathbf{z},\boldsymbol{\zeta}^*)$ it suffices to
compute $\Box_q^{\,l}P^{\mathbf{k}}(\mathbf{z},\boldsymbol{\zeta}^*)$.
We observe that
\begin{equation}\label{shift-prop}
\Box_q\left(\mathbb{C}[M_n]^{\mathbf{k}}_q\right)=\begin{cases}
\mathbb{C}[M_n]^{\mathbf{k}-1}_q, & k_n\geq 1\\ \{0\}, & \mathrm{otherwise}
\end{cases}.
\end{equation}
(We use the notation $\mathbf{k}-l=(k_1-l,k_2-l,\ldots,k_n-l)$.) Indeed, the
operator $\Box_q$ is, in particular, a morphism of
$U_q(\mathfrak{sl}_{n}\times\mathfrak{sl}_{n})$-modules, and all
$U_q(\mathfrak{sl}_{n}\times\mathfrak{sl}_{n})$-modules, isomorphic to
$\mathbb{C}[M_n]^{\mathbf{k}}_q$, have the form
$\mathbb{C}[M_n]^{\mathbf{k}+m}_q$ for some (positive or negative) $m$. The
left hand side is then a subspace of the right hand side by (\ref{knbox}). On the other hand, if
$k_n\geq 1$ and
$\Box_q\left(\mathbb{C}[M_n]^{\mathbf{k}}_q\right)\varsubsetneqq
\mathbb{C}[M_n]^{\mathbf{k}-1}_q$ then
$\Box_q\left(\mathbb{C}[M_n]^{\mathbf{k}}_q\right)=\{0\}$ since
$\mathbb{C}[M_n]^{\mathbf{k}}_q$ and $\mathbb{C}[M_n]^{\mathbf{k}-1}_q$ are
simple isomorphic $U_q(\mathfrak{sl}_{n}\times\mathfrak{sl}_{n})$-modules.
This, however, contradicts positive definiteness of the $q$-Fock
inner product: we have $\Box_q(\det\nolimits_q(\mathbf{z})\cdot f)=0$ for
arbitrary element $f\in\mathbb{C}[M_n]^{\mathbf{k}-1}_q$ and so
$$0=\left(\Box_q(\det\nolimits_q(\mathbf{z})\cdot f),f\right)_{F}=
\left(\det\nolimits_q(\mathbf{z})\cdot f,\det\nolimits_q(\mathbf{z})\cdot
f\right)_{F}.$$

The
equality (\ref{shift-prop}), together with (\ref{detdet}) and (\ref{relinv}), implies
that for certain constant $c(\mathbf{k},l)$
$$\Box_q^{\,l}P^{\mathbf{k}}(\mathbf{z},\boldsymbol{\zeta}^*)=c(\mathbf{k},l)
{P}^{\mathbf{k}-l}(\mathbf{z},\boldsymbol{\zeta}^*)
\mathrm{det}_q(\boldsymbol{\zeta})^{*l}.$$ Indeed, the $q$-integral
operators with the kernels
$\Box_q^{\,l}P^{\mathbf{k}}(\mathbf{z},\boldsymbol{\zeta}^*)$ and
${P}^{\mathbf{k}-l}(\mathbf{z},\boldsymbol{\zeta}^*)
\mathrm{det}_q(\boldsymbol{\zeta})^{*l}$ belong to
$$\mathrm{Hom}_{U_q(\mathfrak{sl}_{n}\times\mathfrak{sl}_{n})}
\left(\mathbb{C}[M_n]^{\mathbf{k}}_q,
\mathbb{C}[M_n]^{\mathbf{k}-l}_q\right),$$ and thus differ by a constant
(the latter space is one-dimensional since $\mathbb{C}[M_n]^{\mathbf{k}}_q$
and $\mathbb{C}[M_n]^{\mathbf{k}-l}_q$ are isomorphic irreducible
$U_q(\mathfrak{sl}_{n}\times\mathfrak{sl}_{n})$-modules). Our immediate aim
is to compute $c(\mathbf{k},l)$.

Given an inner product $(\,\cdot\,,\,\cdot\,)$ in $\mathbb{C}[M_n]_q$ and a
kernel ${P(\mathbf{z},\boldsymbol{\zeta}^*)=\sum_j p'_j\otimes p''_j\in
\mathrm{Pol}(M_n \times S(\mathscr{D}))_q}$ we shall write
$(f(\mathbf{z}),P(\mathbf{z},\boldsymbol{\zeta}^*))$ instead of
$\sum_j(f,p'_j)\cdot (p''_j)^*$.

Recall the notation $(\,\cdot\,,\,\cdot\,)_{S(\mathscr{D})}$ for the inner
product defined in subsection \ref{FFIP} via the invariant integral on the
$q$-Shilov boundary. It is not difficult to observe that the reproducing
property of the kernel $P^{\mathbf{k}}(\mathbf{z},\boldsymbol{\zeta}^*)$ is
equivalent to
$$
(f(\mathbf{z}),P^{\mathbf{k}}(\mathbf{z},\boldsymbol{\zeta}^*))
_{S(\mathscr{D})}=f,\quad \forall f\in \mathbb{C}[M_n]^{\mathbf{k}}_q.$$

We have
$$\Box_q^{\,l}P^{\mathbf{k}}(\mathbf{z},\boldsymbol{\zeta}^*)=c(\mathbf{k},l)
{P}^{\mathbf{k}-l}(\mathbf{z},\boldsymbol{\zeta}^*)\mathrm{det}_q
(\boldsymbol{\zeta})^{*l}$$ or
$$\left(f(\mathbf{z}),\Box_q^{\,l}P^{\mathbf{k}}(\mathbf{z},
\boldsymbol{\zeta}^*)\right)_{S(\mathscr{D})}=
c(\mathbf{k},l)\cdot\left(f(\mathbf{z}),
{P}^{\mathbf{k}-l}(\mathbf{z},\boldsymbol{\zeta}^*)\mathrm{det}_q
(\boldsymbol{\zeta})^{*l}\right)_{S(\mathscr{D})},$$ for any $
f\in\mathbb{C}[M_n]^{\mathbf{k}-l}_q$. By Theorem \ref{qfk}
$$\frac{(1-q^2)^{k_1+k_2+\ldots+k_n-2ln}}{\prod_{i=1}^n
(q^{2n+2-2i};q^2)_{k_i-l}}
\left(f(\mathbf{z}),\Box_q^{\,l}P^{\mathbf{k}}(\mathbf{z},
\boldsymbol{\zeta}^*)\right)_{F}=
c(\mathbf{k},l)\cdot
f(\boldsymbol{\zeta})\cdot\mathrm{det}_q(\boldsymbol{\zeta})^{l},$$ and, due
to the main property of the $q$-Fock product,
$$\frac{(1-q^2)^{k_1+k_2+\ldots+k_n-2ln}}{\prod_{i=1}^n(q^{2n+2-2i};q^2)
_{k_i-l}}
\left(f(\mathbf{z})\mathrm{det}_q(\mathbf{z})^l,P^{\mathbf{k}}(\mathbf{z},
\boldsymbol{\zeta}^*)\right)_{F}=
c(\mathbf{k},l)\cdot
f(\boldsymbol{\zeta})\cdot\mathrm{det}_q(\boldsymbol{\zeta})^{l}.$$ Apply
Theorem \ref{qfk} once again:
$$\frac{\prod_{i=1}^n(q^{2n+2-2i};q^2)_{k_i}}{(1-q^2)^{2ln}\cdot
\prod_{i=1}^n(q^{2n+2-2i};q^2)_{k_i-l}}
\left(f(\mathbf{z})\mathrm{det}_q(\mathbf{z})^l,P^{\mathbf{k}}(\mathbf{z},
\boldsymbol{\zeta}^*)\right)_{S(\mathscr{D})}=
c(\mathbf{k},l)\cdot
f(\boldsymbol{\zeta})\cdot\mathrm{det}_q(\boldsymbol{\zeta})^{l}.$$
The reproducing property of the kernel
$P^{\mathbf{k}}(\mathbf{z},\boldsymbol{\zeta}^*)$ implies that
$$\frac{\prod_{i=1}^n(q^{2n+2-2i};q^2)_{k_i}}{(1-q^2)^{2ln}\cdot
\prod_{i=1}^n(q^{2n+2-2i};q^2)_{k_i-l}} \cdot
f(\boldsymbol{\zeta})\cdot\mathrm{det}_q(\boldsymbol{\zeta})^{l}=
c(\mathbf{k},l)\cdot
f(\boldsymbol{\zeta})\cdot\mathrm{det}_q(\boldsymbol{\zeta})^{l},$$
consequently,
$$c(\mathbf{k},l)=\frac{\prod_{i=1}^n(q^{2n+2-2i};q^2)_{k_i}}
{(1-q^2)^{2ln}\cdot
\prod_{i=1}^n(q^{2n+2-2i};q^2)_{k_i-l}}.$$
We then get
$$\Box_q^{\,l}K_n(\mathbf{z},\boldsymbol{\zeta}^*)=\sum_{\mathbf{k}\,
:\, k_n\geq l}
\frac{\prod_{i=1}^n(q^{2n+2-2i};q^2)_{k_i}}{(1-q^2)^{2ln}\cdot
\prod_{i=1}^n(q^{2n+2-2i};q^2)_{k_i-l}}\cdot
{P}^{\mathbf{k}-l}(\mathbf{z},\boldsymbol{\zeta}^*)
\mathrm{det}_q(\boldsymbol{\zeta})^{*l}=$$
$$=\sum_{\mathbf{k}}\frac{\prod_{i=1}^n(q^{2n+2-2i};q^2)_{k_i+l}}
{(1-q^2)^{2ln}\cdot\prod_{i=1}^n(q^{2n+2-2i};q^2)_{k_i}}\cdot
{P}^{\mathbf{k}}(\mathbf{z},\boldsymbol{\zeta}^*)
\mathrm{det}_q(\boldsymbol{\zeta})^{*l}.$$ Finally, Proposition \ref{coeffi}
implies
$$\Box_q^{\,l}K_n(\mathbf{z},\boldsymbol{\zeta}^*)=
\frac{\prod_{i=1}^n(q^{2n+2-2i};q^2)_{l}}{(1-q^2)^{2ln}}\cdot
K_{n+l}(\mathbf{z},\boldsymbol{\zeta}^*)\cdot
\mathrm{det}_q(\boldsymbol{\zeta})^{*l}.$$ We have thus obtained

\medskip

\begin{proposition}\label{theorem3} For any element $f \in \mathbb{C}[M_n]_q$
one has
$$\Box_q^{\,l}f(\mathbf{z})=\frac{\prod_{i=1}^n(q^{2n+2-2i};q^2)_{l}}
{(1-q^2)^{2ln}}\cdot\int
\limits_{S(\mathscr{D})_q}K_{n+l}(\mathbf{z},\boldsymbol{\zeta}^*)
\mathrm{det}_q(\boldsymbol{\zeta})^{*l}f(\boldsymbol{\zeta})d
\nu(\boldsymbol{\zeta}).$$
\end{proposition}

Proposition \ref{theorem3} reduces the statement of Theorem \ref{Bol} to the
following proposition.

\medskip

\begin{proposition}\label{last} The integral operator
$$ f(\mathbf{z})\mapsto\int
\limits_{S(\mathscr{D})_q}K_{n+l}(\mathbf{z},\boldsymbol{\zeta}^*)
\mathrm{det}_q(\boldsymbol{\zeta})^{*l}f(\boldsymbol{\zeta})d
\nu(\boldsymbol{\zeta})$$ intertwines the $U_q\mathfrak{sl}_{2n}$-actions
$\pi_{n-l}$ and $\pi_{n+l}$.
\end{proposition}
\noindent{\bf Proof.} We shall use a quantum version of the description
(\ref{altdes}) of the twisted action $\pi_\lambda$.

Let us extend the algebra $\mathbb{C}[M_n]_q$ by adding one more generator
$t$ (an analog of the distinguished Pl\"ucker coordinate $t$ in
(\ref{altdes})) such that
$$tz_a^\alpha=q^{-1}z_a^\alpha t, \quad a,\alpha=1,2,\ldots,n.$$ The
localization of
the resulting algebra with respect to the multiplicative system
$t^{\mathbb{N}}$ will be denoted by $\mathbb{C}[M_n]_{q,t}$. It was noted in
\cite{SSV1} that there exists a unique extension of the $U_q
\mathfrak{sl}_{2n}$-module algebra structure in $\mathbb{C}[M_n]_q$ to the
one in $\mathbb{C}[M_n]_{q,t}$ such that
\begin{equation}\label{actio}
E_jt=F_jt=(K_j^{\pm 1}-1)t=0 \quad (j \ne n), \quad
F_nt=(K_n^{\pm1}-q^{\mp1})t=0,\quad E_nt=q^{-1/2}tz_n^n.
\end{equation}
It is clear that the subspace $\mathbb{C}[M_n]_{q}\cdot
t^{-\lambda}\subset\mathbb{C}[M_n]_{q,t}$ is
$U_q\mathfrak{sl}_{2n}$-invariant for any $\lambda\in\mathbb{Z}$. The
following is an equivalent definition of the $U_q\mathfrak{sl}_{2n}$-action
$\pi_\lambda$:
$$(\pi_\lambda(\xi)f)\cdot t^{-\lambda}=\xi(f\cdot t^{-\lambda}), \quad
\xi\in
U_q\mathfrak{sl}_{2n},\quad f\in \mathbb{C}[M_n]_{q}.$$ In other words, the
linear map
$$\mathbb{C}[M_n]_{q}\to\mathbb{C}[M_n]_{q}\cdot
t^{-\lambda},\qquad f\mapsto f\cdot t^{-\lambda}$$ intertwines the
$U_q\mathfrak{sl}_{2n}$-action $\pi_\lambda$ and the natural
$U_q\mathfrak{sl}_{2n}$-action in $\mathbb{C}[M_n]_{q}\cdot t^{-\lambda}$.

We need also certain extension of the algebra
$\mathrm{Pol}(S(\mathscr{D}))_q$. Let us add to
$\mathrm{Pol}(S(\mathscr{D}))_q$ two generators $t, t^*$ such that $$
tt^*=t^*t,\quad tz_a^\alpha=q^{-1}z_a^\alpha t,\quad
t^*z_a^\alpha=q^{-1}z_a^\alpha t^*,\quad a,\alpha=1,2,\ldots,n.$$ Denote
this new algebra by $\mathrm{Pol}(\widehat{S}(\mathscr{D}))_q$ and its
localization with respect to the multiplicative system $(tt^*)^{\mathbb{N}}$
by $\mathrm{Pol}(\widehat{S}(\mathscr{D}))_{q,x}$. The involution in
$\mathrm{Pol}(S(\mathscr{D}))_q$ can be extended to an involution in
$\mathrm{Pol}(\widehat{S}(\mathscr{D}))_{q,x}$ by setting $*:t\mapsto t^*$.
It is proved in \cite[Section 2]{CS} that there exists a unique structure
of $U_q \mathfrak{su}_{n,n}$-module algebra in
$\mathrm{Pol}(\widehat{S}(\mathscr{D}))_{q,x}$ which coincides with the
original one on
$\mathrm{Pol}({S}(\mathscr{D}))_q\subset\mathrm{Pol}(\widehat{S}
(\mathscr{D}))_{q,x}$ and satisfies (\ref{actio}). Following \cite[Section
3]{CS}, we equip $\mathrm{Pol}(\widehat{S}(\mathscr{D}))_{q,x}$ with a
bigrading: $$ \deg t=(0,1),\quad \deg t^*=(1,0),\quad
\deg(z_a^\alpha)=\deg(z_a^\alpha)^*=(0,0),\quad a,\alpha=1,2,\ldots,n. $$
Obviously, the homogeneous components $$
\mathrm{Pol}(\widehat{S}(\mathscr{D}))_{q,x}^{(i,j)}=\{f \in
\mathrm{Pol}(\widehat{S}(\mathscr{D}))_{q,x}|\:\deg f=(i,j)\}=t^{*i}\cdot
\mathrm{Pol}(S(\mathscr{D}))_q \cdot t^j $$ are submodules of the $U_q
\mathfrak{sl}_{2n}$-module $\mathrm{Pol}(\widehat{S}(\mathscr{D}))_{q,x}$.

Proposition \ref{last} is an immediate consequence of the following
statement

\medskip

\begin{lemma}The linear operator from
$\mathrm{Pol}(\widehat{S}(\mathscr{D}))_{q,x}^{(0,l-n)}$ to
$\mathbb{C}[M_n]_{q}\cdot t^{-l-n}$ given by $$f\cdot t^{l-n} \mapsto
\left(\int
\limits_{S(\mathscr{D})_q}K_{n+l}(\mathbf{z},\boldsymbol{\zeta}^*)
\mathrm{det}_q(\boldsymbol{\zeta})^{*l}f(\boldsymbol{\zeta})d
\nu(\boldsymbol{\zeta})\right)\cdot t^{-l-n}$$ is a morphism of
$U_q\mathfrak{sl}_{2n}$-modules.
\end{lemma}
\noindent{\bf Proof of the lemma.} The proof may be easily reduced to the
following three statements:

\noindent (i) The linear map
\begin{equation}\label{map1}
\mathrm{Pol}(\widehat{S}(\mathscr{D}))_{q,x}^{(0,l-n)}\to
\mathrm{Pol}(\widehat{S}(\mathscr{D}))_{q,x}^{(l,-n)},\quad
f(\boldsymbol{\zeta})\cdot t^{l-n}\mapsto
\mathrm{det}_q(\boldsymbol{\zeta})^{*l}f(\boldsymbol{\zeta})\cdot
t^{*l}t^{-n}
\end{equation}
is a morphism of $U_q\mathfrak{sl}_{2n}$-modules. This statement follows from the results of Sections 2,3 in \cite{CS}.

\noindent (ii) Let
$K_{n+l}(\mathbf{z},\boldsymbol{\zeta}^*)=\sum_jk'_j(\mathbf{z})\otimes
k''_j(\boldsymbol{\zeta}^*)$. Then the element $$\sum_jk'_j(\mathbf{z})\cdot
t^{-l-n}\otimes t^{*(-l-n)}\cdot k''_j(\boldsymbol{\zeta}^*)\in
(\mathbb{C}[M_n]_{q}\cdot
t^{-l-n})\widehat{\otimes}\mathrm{Pol}(\widehat{S}(\mathscr{D}))_{q,x}^
{(-l-n,0)}$$ is a $U_q\mathfrak{sl}_{2n}$-invariant (here the symbol
$\widehat{\otimes}$ has the same meaning as the one in the equality
$\mathbb{C}[M_n]_q\widehat{\otimes}
\mathrm{Pol}(S(\mathscr{D}))_q=\mathrm{Fun}(M_n \times S(\mathscr{D}))_q$).
The statement is a consequence of results of Section 8 in \cite{SSV2}. This,
together with statement (i), implies that the map
$$\mathrm{Pol}(\widehat{S}(\mathscr{D}))_{q,x}^{(0,l-n)}\to
(\mathbb{C}[M_n]_{q}\cdot
t^{-l-n})\widehat{\otimes}\mathrm{Pol}(\widehat{S}(\mathscr{D}))_{q,x}
^{(-n,-n)},$$
\begin{equation}\label{map2}
f(\boldsymbol{\zeta})\cdot t^{l-n}\mapsto
\sum_jk'_j(\mathbf{z})\cdot t^{-l-n}\otimes t^{*(-l-n)}\cdot
k''_j(\boldsymbol{\zeta}^*)\mathrm{det}_q(\boldsymbol{\zeta})^{*l}
f(\boldsymbol{\zeta})\cdot
t^{*l}t^{-n}
\end{equation}
is a morphism of $U_q\mathfrak{sl}_{2n}$-modules.

\noindent (iii) The linear functional
$$\mathrm{Pol}(\widehat{S}(\mathscr{D}))_{q,x}^{(-n,-n)}\to\mathbb{C},\quad
t^{*(-n)}\cdot f\cdot t^{-n} \mapsto \int
\limits_{S(\mathscr{D})_q}f(\boldsymbol{\zeta})d \nu(\boldsymbol{\zeta})$$
is a $U_q\mathfrak{sl}_{2n}$-invariant integral. This is proved in Section 3
of \cite{CS}. As a consequence of this statement and statements (i), (ii) we
get: the linear operator from
$\mathrm{Pol}(\widehat{S}(\mathscr{D}))_{q,x}^{(0,l-n)}$ to
$\mathbb{C}[M_n]_{q}\cdot t^{-l-n}$ given by
\begin{equation}\label{map3}
f(\boldsymbol{\zeta})\cdot t^{l-n}\mapsto \sum_jk'_j(\mathbf{z})\cdot
t^{-l-n}\cdot \int \limits_{S(\mathscr{D})_q}\Theta_l\left(
k''_j(\boldsymbol{\zeta}^*)\mathrm{det}_q(\boldsymbol{\zeta})^{*l}f
(\boldsymbol{\zeta})\right)d
\nu(\boldsymbol{\zeta})
\end{equation}
is a morphism of $U_q\mathfrak{sl}_{2n}$-modules (here $\Theta_l$ means the
automorphism of the algebra $\mathrm{Pol}({S}(\mathscr{D}))_q$ given by
$f\mapsto t^{*(-l)}\cdot f\cdot t^{*l}$).

Lemma 4.10 follows from the latter statement and the equality
$$\int \limits_{S(\mathscr{D})_q}\Theta_l(f(\boldsymbol{\zeta}))d
\nu(\boldsymbol{\zeta})=\int
\limits_{S(\mathscr{D})_q}f(\boldsymbol{\zeta})d
\nu(\boldsymbol{\zeta}),\quad\forall f$$
which is due to the simple observation that the functional $\int
\limits_{S(\mathscr{D})_q}f(\boldsymbol{\zeta})d
\nu(\boldsymbol{\zeta})$ 'picks up' the constant term of $f$ and the constant terms of $\Theta_l(f)$ and $f$ are the same.
\hfill$\blacksquare$

Proposition \ref{last}, consequently Theorem \ref{Bol}, is
now proved.\hfill$\blacksquare$

\bigskip

\section{Holomorphic discrete series for $U_q\mathfrak{su}_{n,n}$} In this
last section we study the holomorphic discrete series
representations for $U_q\mathfrak{su}_{n,n}$ and study their analytic
continuation.

After giving a definition of the holomorphic discrete series
for $U_q\mathfrak{su}_{n,n}$, we prove an analog of a classical known result by J.~Faraut
and A.~Koranyi which allows one to express the inner product in a module of
the holomorphic discrete series via the Fock inner product. We use
the result to prove unitarizability of the modules $\mathcal{P}_\lambda$
with $\lambda>n-1$;
the discrete series parameters are  $\lambda>2n-1$.
We apply then the covariance property, proved earlier, to studying certain
quotients of the modules $\mathcal{P}_\lambda$.

\medskip

\subsection{Definition of the holomorphic discrete series}

We start by recalling the definition of the holomorphic discrete series
for $\widetilde{SU}_{n,n}$. Fix $\lambda>2n-1$ and consider the Hilbert
space of holomorphic functions on the unit matrix ball $\mathscr{D}$ which
are square integrable with respect to the measure $\det
\nolimits(1-\mathbf{z}\mathbf{z}^*)^{\lambda-2n}d\mathbf{z}$ (here
$d\mathbf{z}$ is the normalized Lebesgue measure:
$\int_{\mathscr{D}}d\mathbf{z}=1$). It is known \cite{Arazy} that the
operators (\ref{weightaction}) are unitary on that Hilbert space. The
corresponding representation of $\widetilde{SU}_{n,n}$ is said to be a
representation of the holomorphic discrete series.

Now let us turn to the quantum setting. Suppose $A_0=(A,*)$ is a Hopf
$*$-algebra. An $A$-module $V$ is said to be a unitarizable $A_0$-module if
there exists an inner product $(\,,\,)$ on $V$ such that for all $v_1,v_2\in
V$ and any $\xi\in A$
$$(\xi v_1,v_2)=(v_1,\xi^*v_2)$$ (that is, $V$ possesses an $A_0$-invariant
inner product, see subsection \ref{FFIP}).

Unitarizable $U_q\mathfrak{su}_{n,n}$-modules substitute unitary
representations of $SU_{n,n}$ (or $\widetilde{SU}_{n,n}$) in the quantum
setting. The following statement was proved in \cite{SSV2}, Corollary 6.5.

\medskip

\begin{proposition}\label{unitary}
For $\lambda>2n-1$ there exists a unique inner product $(\,,\,)_\lambda$ in
$\mathcal{P}_\lambda$ such that for all $f_1,f_2\in \mathcal{P}_\lambda$ and
$\xi\in U_q\mathfrak{sl}_{2n}$
$$
(\pi_\lambda(\xi)f_1,f_2)_\lambda=(f_1,\pi_\lambda(\xi^*)f_2)_\lambda,
$$
with the normalization $(1, 1)_\lambda=1$.
\end{proposition}

\medskip

Clearly, the unitarizable $U_q\mathfrak{su}_{n,n}$-modules
$\mathcal{P}_\lambda$, $\lambda>2n-1$, are $q$-analogs of unitary
representation of the holomorphic discrete series for
$\widetilde{SU}_{n,n}$.

The inner product $(\,,\,)_\lambda$ may be described explicitly as follows.
Let us use the notation from the proof of Proposition \ref{qfk}. Consider
the element
\begin{equation}\label{elem}
\hat{K}_{\lambda}=\frac{\prod_{j=0}^\infty
\left(1+\sum_{k=1}^m(-q^{2(\lambda+j)})^k \hat{\chi}_k
\right)}{\prod_{j=0}^\infty \left(1+\sum_{k=1}^m(-q^{2j})^k \hat{\chi}_k
\right)}\in\mathbb{C}[[M_n \times M_n]]_q.
\end{equation}
Then the inner product $(\,,\,)_\lambda$ is the one associated with the
above element, i.e.
\begin{equation}\label{elem-1}
f=\sum_jk'_j\cdot(f\,,\,k''_j)_\lambda
\end{equation}
provided $\hat{K}_{\lambda}=\sum_jk'_j\otimes k''_j$ (see Theorem 9.1 in
\cite{SSV2}).

\subsection{A $q$-analog of a result by J.~Faraut and
A.~Koranyi}\label{FaKo}

In this subsection, we present a $q$-analog of Corollary 3.7 in \cite{FK}
where (in the classical setting) the Fock inner product and the
inner products $(\,,\,)_\lambda$ are compared. We then apply the result to
the problem of analytic continuation of the holomorphic discrete series for
$U_q\mathfrak{su}_{n,n}$.

Using the arguments preceding Proposition \ref{qfk} we deduce that the inner
product $(\,,\,)_\lambda$ and the $q$-Fock inner product on a
particular simple
$U_q\mathfrak{s}(\mathfrak{gl}_n\times\mathfrak{gl}_n)$-submodule
$\mathbb{C}[M_n]^{\mathbf{k}}_q\subset\mathbb{C}[M_n]_q$ are proportional.
The proportionality constant is given by the following formula.

\medskip

\begin{proposition}\label{qfklamb} Let $\lambda >2n-1$. Then
\begin{equation}\label{lamb}
(f_1,f_2)_{F}=\frac{\prod_{i=1}^n(q^{2\lambda+2-2i};q^2)_{k_i}}
{(1-q^2)^{k_1+k_2+\ldots+k_n}}\cdot(f_1,f_2)_\lambda, \quad
f_1,f_2\in\mathbb{C}[M_n]^{\mathbf{k}}_q.
\end{equation}
\end{proposition}
\noindent{\bf Proof.} Consider the reproducing kernel
$$
{K}_{\lambda}=\frac{\prod_{j=0}^\infty
\left(1+\sum_{k=1}^m(-q^{2(\lambda+j)})^k {\chi}_k
\right)}{\prod_{j=0}^\infty \left(1+\sum_{k=1}^m(-q^{2j})^k {\chi}_k
\right)}
$$
associated to the element (\ref{elem}). The image of this kernel under the
isomorphism (\ref{isom}) is given by
$\prod_{i=1}^{n}\frac{(q^{2\lambda}x_i;q^2)_\infty}{(x_i;q^2)_\infty}$.
By repeating the computation from subsection \ref{kernels}, one gets
$$(f_1,f_2)_{S(\mathscr{D})}=\prod_{i=1}^n\frac{(q^{2\lambda+2-2i};q^2)_{k_i}}
{(q^{2n+2-2i};q^2)_{k_i}}\cdot(f_1,f_2)_\lambda, \quad
f_1,f_2\in\mathbb{C}[M_n]^{\mathbf{k}}_q.$$ What remains is to apply
Proposition \ref{qfk}. \hfill$\blacksquare$

\medskip

The result of the above proposition has an important application to the
problem of analytic continuation of the holomorphic discrete series. The
point is that (\ref{lamb}) allows one to define the sesquilinear form
$(\,,\,)_\lambda$ on $\mathbb{C}[M_n]_q$ for any $\lambda\in\mathbb{R}$ for
which all the multipliers $\prod_{i=1}^n(q^{2\lambda+2-2i};q^2)_{k_i}$ are
non-zero. It is not difficult to prove that the resulting form is still
$U_q\mathfrak{su}_{n,n}$-invariant with respect to the corresponding twisted
action. Indeed, the invariance is equivalent to the infinitely many
equalities of the form
$$(\pi_\lambda(\xi)f_1,f_2)_\lambda=(f_1,\pi_\lambda(\xi^*)f_2)_\lambda,
\quad \xi\in U_q\mathfrak{sl}_{2n},\, f_1,f_2\in\mathbb{C}[M_n]_q.$$ After
simple transformations, each equality becomes an equality of two Laurent
polynomials in $q^{\lambda}$ which is known to hold for $\lambda>2n-1$ by
Proposition \ref{unitary}, and, thus, for any $\lambda$. It is naturally to
pose the problem of finding those $\lambda$'s for which the
corresponding sesquilinear form is positive definite, i.e. the corresponding
$U_q\mathfrak{su}_{n,n}$-modules are unitarizable. In the classical setting,
such $\lambda$'s are said to belong to the continuous part of the Wallach
set \cite{Arazy}.

The proposition implies positive definiteness of the inner product
$(f_1,f_2)_\lambda$ for any $\lambda>n-1$:

\medskip

\begin{corollary}\label{cor}
The $U_q\mathfrak{su}_{n,n}$-modules $\mathcal{P}_\lambda$ are unitarizable
for $\lambda>n-1$.
\end{corollary}

\medskip
\noindent For $n=2$, this statement was obtained in \cite{GeomReal} (see
Proposition 6.1).

\subsection{Some consequences of the covariance property}

In the previous subsection, we were able to deduce some irreducibility
and unitarity property of
$U_q\mathfrak{sl}_{2n}$-modules $\mathcal{P}_\lambda$ from the results obtained earlier. We will not pursue all the details
here. It is immediate that
$\mathcal{P}_\lambda$
is reducible for
$\lambda=n-1, n-2,\ldots$
for the following obvious reason: by the covariance property,
$\mathcal{P}_{n-l}^{(0)}=\mathrm{Ker}\,\Box_q^{\,l}$ is a submodule in
$\mathcal{P}_{n-l}$. A related
 application of the covariance property is the following
\medskip

\begin{proposition}
For any $l\in\mathbb{N}$ $\mathcal{P}_{n-l}/\mathcal{P}_{n-l}^{(0)}$ is a
unitarizable $U_q\mathfrak{su}_{n,n}$-module isomorphic to
$\mathcal{P}_{n+l}$.
\end{proposition}

\medskip
\noindent {\bf Proof.} Unitarizability follows from the covariance property,
Corollary \ref{cor}, and injectivity of the induced morphism
$\mathcal{P}_{n-l}/\mathcal{P}_{n-l}^{(0)}\to\mathcal{P}_{n+l}$ of
$U_q\mathfrak{sl}_{2n}$-modules. Actually, this latter morphism is an
isomorphism. To prove this, it suffices to show that the operator
$$\Box_q^{\,l}:\mathbb{C}[M_n]_q\to\mathbb{C}[M_n]_q$$ is surjective.
In turn, it suffices to prove the latter statement for $l=1$, i.e. to show
that $\Box_q$ is surjective. But this follows
from (\ref{shift-prop}).
 \hfill$\blacksquare$

\end{document}